\documentclass{amsart}
\usepackage[latin1]{inputenc}
\usepackage{amssymb}
\usepackage{amsmath}
\usepackage{enumerate}
\usepackage{epsfig}
\usepackage{subfigure}
\usepackage[all]{xy}
\usepackage{longtable}
\usepackage{tikz}

\newtheorem{theorem}{Theorem}[section]
\newtheorem{lemma}[theorem]{Lemma}
\newtheorem{proposition}[theorem]{Proposition}
\newtheorem{corollary}[theorem]{Corollary}

\theoremstyle{definition}
\newtheorem{definition}[theorem]{Definition}

\theoremstyle{remark}
\newtheorem{rem}[theorem]{Remark}
\theoremstyle{remark}

\numberwithin{equation}{section}

\newcommand{\rk}{{\rm rank}}

\newcommand{\ra}{\rightarrow}

\newcommand{\C}{\mathbb{C}}
\newcommand{\R}{\mathbb{R}}
\newcommand{\Z}{\mathbb{Z}}
\newcommand{\Q}{\mathbb{Q}}

\newcommand{\N}{\mathbb{N}}

\newcommand{\Aut}{{\operatorname{Aut}}}

\makeatletter
\def\blfootnote{\xdef\@thefnmark{}\@footnotetext}

\begin{document}

\keywords{K3 surfaces, Kummer surfaces, Kummer type lattice, Quotient surfaces\\
The author is partially supported by PRIN 2010--2011 ``Geometria
delle variet\`a algebriche" and FIRB 2012 ``Moduli Spaces and
their Applications".}

\title[On K3 surface quotients of K3 or Abelian surfaces]{On K3 surface quotients of K3 or Abelian surfaces}
\author{Alice Garbagnati}
\address{Alice Garbagnati, Dipartimento di Matematica, Universit\`a di Milano,
  via Saldini 50, I-20133 Milano, Italia}
\email{alice.garbagnati@unimi.it}

\begin{abstract}
The aim of this paper is to prove that a K3 surface is the minimal model of the quotient of an Abelian surface by a group $G$ (respectively of a K3 surface by an Abelian group $G$) if and only if a certain lattice is primitively embedded in its N\'eron--Severi group. This allows one to describe the coarse moduli space of the K3 surfaces which are (rationally) $G$-covered by Abelian or K3 surfaces (in the latter case $G$ is an Abelian group).
If either $G$ has order 2 or $G$ is cyclic and acts on an Abelian surface, this result was already known, so we extend it to the other cases.

Moreover, we prove that a K3 surface $X_G$ is the minimal model of the quotient of an Abelian surface by a group $G$ if and only if a certain configuration of rational curves is present on $X_G$. Again this result was known only in some special cases, in particular if $G$ has order 2 or 3.
\end{abstract}

\maketitle

\section{Introduction}
Thanks to the Torelli theorem for K3 surfaces and to the theory of the lattice polarized K3 surfaces, in order to describe the moduli space of K3 surfaces having a certain geometric property it is useful to translate this geometric property in terms of embeddings of certain lattices. In this paper we analyze the geometric property "a K3 surface is the minimal model of the quotient of an Abelian or a K3 surface by a finite group". Under certain conditions we are able to translate this property to a lattice theoretic property and thus to describe the coarse moduli space of the K3 surfaces which are (rationally) covered by Abelian surfaces or by K3 surfaces.
This generalizes several previous results by Nikulin, \cite{NikKummer}, Bertin, \cite{B}, and by Sarti and the author, \cite{GSinv}.

The first and crucial example is given by the Kummer surfaces: a Kummer surface is a K3 surface obtained as minimal resolution of the quotient $A/\iota$, where $A$ is an Abelian surface and $\iota$ is an involution on $A$. In \cite{NikKummer}, Nikulin proved that a K3 surface is a Kummer surface if and only if at least one of the two following equivalent conditions holds:\\
$i)$ a certain lattice, called Kummer lattice, is primitively embedded in the N\'eron--Severi group of the K3 surface;\\
$ii)$ there are 16 disjoint smooth rational curves on the K3 surface.\\
The first condition is more related with the lattice theory and allows one to describe the coarse moduli space of the K3 surfaces which are Kummer surfaces. The second one is clearly more related with the geometry of the surface.

In a more general setting, we are considering the following situation: $Y$ is either an Abelian or a K3 surface, $G$ is a finite group of automorphisms of $Y$ and the minimal model of $Y/G$ is a K3 surface $X$. In this case we say that $X$ is (rationally) $G$-covered by $Y$. In view of the results by Nikulin on Kummer surfaces, it is quite natural to pose the following two questions:\\
 
{\bf Question A)} Is the property "a K3 surface $X$ is (rationally) $G$-covered by a surface $Y$" equivalent to the condition "there is a certain lattice (depending on $G$) which is primitively embedded in $NS(X)$"?\\
We observe that a positive answer to this question immediately provides a description of the coarse moduli space of the K3 surfaces (rationally) covered by Abelian or K3 surfaces.\\

{\bf Question B)} Is the property "a K3 surface $X$ is (rationally) $G$-covered by a surface $Y$" equivalent to the condition "there is a certain configuration of rational curves on $X$"?\\

The main results of this paper are to give a positive answer to\begin{itemize}\item  question A) in the cases $Y$ is an Abelian surface (see Theorem \ref{theorem: X is A/G iff KG is in NS(X)}); \item question B) in the case $Y$ is an Abelian surface (see Theorem \ref{theor: it suffices to have the curves on X});
\item question A) in the cases $Y$ is a K3 surface and $G$ is an Abelian group (see Theorem \ref{theorem: S is A/G iff MG is in NS(S)}).\end{itemize}

It is not possible to give in general a positive answer to question B) in the case $Y$ is a K3 surface. For example, it is known that the answer is negative if we assume that $Y$ is a K3 surface and $G=\Z/2\Z$, cf. \cite{GS}.
I do not know if it is possible to extend the positive answer given to the question A) in the case that $Y$ is a K3 surface and $G$ is an Abelian group to the weaker hypothesis that $Y$ is a K3 surface, without assumptions (or with different assumptions) on $G$. \\

The positive answer to the question A) in case $Y$ is an Abelian surface was already known if $G$ is a cyclic group, indeed the classical case of the Kummer surface, i.e. $G=\Z/2\Z$, was considered by Nikulin, \cite{NikKummer}, as we said above, the other cyclic cases are considered in \cite{B}. In Theorem \ref{theorem: X is A/G iff KG is in NS(X)} we do the remaining cases. In order to state and prove this theorem, the first step is to find all the finite groups $G$ acting on an Abelian surface, in such a way that $A/G$ desingularizes to a K3 surface. We assume, without loss of generalities, that $G$ does not contain translations. The list of these groups is classically known, see \cite{F}, and consists of 4 cyclic groups and 3 non cyclic (and non Abelian) groups. One of the non cyclic groups, the quaternion group, can act on two different families of Abelian surfaces and the actions have different sets of points with non trivial stabilizer. So we have to consider 4 actions of non cyclic groups on an Abelian surface. The second step is the identification of the lattice that should characterize the K3 surfaces which are (rationally) $G$-covered by an Abelian surface. This lattice depends on $G$ and the natural candidate (also in view of the previous results by Nikulin and Bertin) is the minimal primitive sublattice of $NS(X)$ which contains all the curves arising from the desingularization of $A/G$. We call the lattices constructed in this way lattices of Kummer type and we denote them by $K_G$. They were already determined if $G$ is a cyclic group and they are computed in the non-cyclic case in Section \ref{subsection: non cyclic quotients of A}. The lattices arising in the non cyclic cases were already considered in \cite{W}, but unfortunately some of the results presented in \cite[Proposition 2.1]{W} are not correct.
The last step is to prove our main result, that is: the primitive embedding of the lattices $K_G$ in the N\'eron--Severi group of a K3 surface $X$ is equivalent to the fact that $X$ is (rationally) $G$-covered by an Abelian surface. We combine a result of \cite{B}, which allows us to give a geometric interpretation of the $(-2)$-classes appearing in the lattice $K_G$, with classical results on cyclic covers between surfaces, in order to reconstruct the surface $A$ starting from $(X,K_G)$.

The positive answer to the question A) in case $Y$ is a K3 surface and $G$ is an Abelian group, is contained in Theorem \ref{theorem: S is A/G iff MG is in NS(S)}. The admissible Abelian groups are listed in \cite[Theorem 4.5]{Niksympl} and there are 14. Only the case $G=\Z/2\Z$ was already known, see \cite{GSinv}. The proof of the result is totally analogous to the one given in case $Y$ is an Abelian surface, with the advantage that the lattices of Kummer type are substituted by other lattices, denoted by $M_G$, which were already computed in \cite[Sections 6 and 7]{Niksympl} for all the admissible groups $G$. The extension to the non Abelian groups $G$ seems more complicated: the lattices $M_G$ are known also in the non Abelian cases, \cite{X}, but  it is not so clear how to reconstruct the surface $Y$ only from the data $(X, M_G)$.

Let us now discuss the more geometric question B). The positive answer to question B) in case $Y$ is an Abelian surface was already known in cases $G=\Z/2\Z$ and $G=\Z/3\Z$. The case of the involution was considered by Nikulin, as we said, in \cite{NikKummer}. The case $G=\Z/3\Z$ is due to Barth, \cite{Barthnine}. The other groups are considered here. The proof of this very geometric result is essentially based on computations in lattice theory. Indeed the idea is to prove that if a K3 surface $X$ admits a certain configuration of curves, then the minimal primitive sublattice of the N\'eron--Severi group containing these curves is in fact $K_G$. We underline that the computations with these lattices are strongly conditioned by that we are considering many curves, which implies that the rank of the lattices that they span is high. This is exactly the hypothesis which fails if we consider the case $Y$ is a K3 surface (and not an Abelian surface). Indeed in this case the result can not be extended (at least without conditions on the Abelian group $G$).

In Section \ref{sec: preliminaries} we recall some known results. In Section \ref{sec: Bertin  lemma and implications} we present Proposition \ref{prop: the roots of L are smooth  irreducible curves} (based on previous results by Bertin) which is fundamental in the proof of our main theorems. In Section \ref{subsection: non cyclic quotients of A} we compute the lattices $K_G$ in case $G$ is not a cyclic group and in Section \ref{subsec: the Kummer type  lattices} we give an exhaustive description of the lattices of Kummer type and of their properties. In Section \ref{subsec:  main on K3 covered by Abelian} we state and prove two of our main results, giving a positive answer to the Questions A and B in case $Y$ is an Abelian surface (see Theorem \ref{theorem: X is A/G iff KG is in NS(X)} and Theorem \ref{theor: it suffices to have the curves on X}).  In Section \ref{K3 rationally covered by Z/3Z} we discuss the relation between K3 surfaces which are (rationally) $\Z/3\Z$-covered by Abelian surfaces, K3 surfaces which are (rationally) $(\Z/3\Z)^2$-covered by K3 surfaces, K3 surfaces which are (rational) $(\Z/3\Z)^2$-covers of K3 surface. This generalizes a similar result on Kummer surfaces, proved in \cite{GS}.\\
In Section \ref{sc: K3 covered  by K3} we concentrate ourselves on K3 surfaces covered by K3 surfaces, giving a positive answer to question A) in this settings. Moreover, we give more precise results on the K3 surfaces which a (rationally) $\Z/3\Z$-covered by K3 surfaces, presenting all the possible N\'eron--Severi groups of a K3 surface with this property and minimal Picard number. This generalizes results proved in \cite{GSinv}.

{\bf Notation}:\\
$\bullet$ $\mathbb{D}_{4n}$ is the dicyclic group of
order $4n$ (called also binary dihedral group), which has the following presentation: $\langle a,b,c| a^n=b^2=c^2=abc\rangle$. (In case $n=2$ it is the quaternion group). In the GAP ID, $\mathbb{D}_8$ is $(8, 4)$ and $\mathbb{D}_{12}$ is $(12,1)$;\\
$\bullet$ $\mathbb{T}$ is the binary tetrahedral group: it has
order 24 and the following presentation $\langle
r,s,t|r^2=s^3=t^3=rst\rangle$. In the GAP ID, $\mathbb{T}$ is $(24, 3)$;\\
$\bullet$ $\mathfrak{A}_{3,3}$ is the generalized dihedral group of the abelian group of order 9: it has order 18 and the following presentation: $\langle
r,s,t|r^2=s^3=t^3=1,\ tr=rt^2,\ sr=rs^2\rangle$. In the GAP ID, $\mathfrak{A}_{3,3} $ is $(18, 4)$ 

\section{Preliminaries}\label{sec: preliminaries}
In this Section we recall some very well known facts and fix the notation.

\subsection{Lattices}\label{subsec: lattice}
\begin{definition} A lattice is a pair $(L,b_L)$, where $L=\Z^n$, $n\in \N$ and $b_L:L\times L\ra\Z$ is a symmetric non degenerate bilinear form taking values in $\Z$. The number $n$ is the rank of $L$. The signature of $(L,b_L)$ is the signature of the $\R$ linear extension of $(L,b_L)$.

A lattice is said to be even if the quadratic form induced by $b_L$ on $L$ takes value in $2\Z$, and not only in $\Z$.

The discriminant group of $L$ is $L^{\vee}/L$, where the dual $L^{\vee}$ can be identified with the set $\{m\in L\otimes \Q\ |\ b_L(m,l)\in \Z \mbox{ for all }l\in\Z\}$ (here we denote by $b_L$ also the $\Q$ linear extension of $b_L$). The discriminant form is the form induced by $b_L$ on the discriminant group. 

The length of a lattice $(L,b_L)$, denoted by $l(L)$, is the minimal number of generators of the discriminant group.

A lattice is said to be unimodular if its discriminant group is trivial, i.e. if its length is zero.\end{definition}

The discriminant group of a lattice is a finite group, free product of cyclic groups. Its order is the determinant of one (and so of any) matrix which represents the form $b_L$ with respect to a certain basis of $L$. This number is called discriminant of the lattice $L$ and is denoted by $d(L)$.

In the following we are interested in the construction of overlattices of finite index of a given lattice. Let $L$ and $M$ be two lattices with the same rank. Let $L\hookrightarrow M$. Then $M$ is generated by the vectors which generate $L$ and by some other vectors, which are non trivial in $M/L$ but which necessarily have an integer intersection with all the vectors in $L$ (otherwise the form on $M$ can not take values in $\Z$). This means that the non trivial  vectors in $M/L$ are non trivial elements of in the discriminant group of $L$. 

If moreover we require that the lattice $M$ is even, then $L$ is automatically even (since it is a sublattice of $M$) and also the non trivial classes in $M/L$ have an even self intersection. So, if we have an even lattice $L$ and we want to construct an even overlattice of finite index, we have to add to the generators of $L$ some non trivial elements in $L^{\vee}/L$ which have an even self intersection.

More in general, every isotropic subgroup of $L^{\vee}/L$ (where a subgroup $H$ is isotropic if the discriminant form restricted to $H$ is trivial), corresponds to an overlattice of finite index of $L$ and viceversa every overlattice of finite index of $L$ corresponds to an isotropic subgroup of $L^{\vee}/L$ , see Section \cite[Section 1]{Nikbilinear}

If $M$ is an overlattice of $L$ of index $r$, then $d(L)/d(M)=r^2$.

\begin{definition} Let $M$ and $L$ be two lattices with $\rk(M)\leq \rk(L)$. Let $\varphi:M\ra L$ be an embedding of $M$ in $L$. We say that $\varphi$ is primitive, or that $M$ is primitively embedded in $L$, if $L/\varphi(M)$ is torsion free. \end{definition}

\begin{proposition}\label{prop: length of lattices}{\rm (see \cite[Proposition 1.6.1]{Nikbilinear})} Let $L$ be a unimodular lattice, $M$ be a primitively sublattice of $L$ and $M^{\perp}$ the orthogonal to $M$ in $L$. The discriminant group of $M$ is isomorphic to the discriminant group of $M^{\perp}$. In particular, since the length of a lattice is at most the rank of the lattice, $l(M)=l(M^{\perp})\leq \min\{\rk(M), \rk(L)-\rk(M^{\perp})\}$.\end{proposition}

\begin{definition} A root of the lattice $(L,b_L)$ is a vector $v\in L$ such that $b_L(v,v)=-2$. The root lattice of a given lattice $L$ is the lattice spanned by the set of all the roots in $L$.

A lattice is called root lattice if it is generated by its roots. In particular a root lattice is negative definite.\end{definition}

\subsection{Covers}\label{sec: covers}

Here we recall a very well known and classical result (see \cite[Chapter I, Section 17]{BHPV}) on covers, which will be essential for our purpose.

Let $Y$ be a connected complex manifold and $B$ an effective divisor on $Y$. Suppose we have a line bundle $\mathcal{L}$ on $Y$ such that \begin{equation}\label{eq: divisible divisor}\mathcal{O}_Y(B)=\mathcal{L}^{\otimes n},\end{equation}
and a section $s\in H^0(Y,\mathcal{O}_Y(B))$ vanishing exactly along $B$. We denote by $L$ the total space of $\mathcal{L}$ and we
let $p:L\ra  Y$ be the bundle projection. If $t\in H^0(L,p^*(\mathcal{L}))$ is the tautological
section, then the zero divisor of $p^*s -t$ defines an analytic subspace $X$ in $L$. The variety $X$ is an $n$-cover of $Y$ branched along $B$ and determined by $\mathcal{L}$ and the cover map is the restriction of $p$ to $X$.
If both $Y$ and $B$ are  smooth and reduced, then $X$ is smooth.

Let us denote by $D\in Pic(Y)$ the divisor associated to the line bundle $\mathcal{L}$. The condition \eqref{eq: divisible  divisor} is equivalent to $B=nD$, i.e. $B/n=D\in Pic(Y)$. For this reason we call $B$ an $n$-divisible divisor (or set) in the Picard group. We often call the curves in the support of $B$ an $n$-divisible set of curves. The previous discussion implies that to each effective divisible divisor one can associate a cyclic cover of the variety.

Let us consider a sort of viceversa: let $\pi:X\ra Y$ an  $n$-cyclic cover between smooth varieties such that the branch locus is smooth and all its components have codimension 1 in $Y$. Then $\pi$ determines a divisor (supported on the branch locus) which is divisible by $n$. 
This applies in particular to a special situation that we will considered in the following: let $X$ and $Y$ be two surfaces. Let $\alpha$ be an automorphism of $X$ of order 2 or 3 which fixes only isolated points. Then it is possible to construct a blow up $\widetilde{X}$ of $X$ such that $\alpha$ induces an automorphism $\widetilde{\alpha}$ of $\widetilde{X}$ whose fixed locus consists of disjoint curves. So $\widetilde{X}/\widetilde{\alpha}$ is a smooth surface that we denote by $Y$ and which is birational to $X/\alpha$. The quotient map $\widetilde{X}\ra \widetilde{X}/\widetilde{\alpha}=Y$ is a $|\alpha|:1$ cover of $Y$ branched along a smooth union of curves. Hence there is a $|\alpha|$-divisible set of curves on $Y$, which is a smooth birational model of $X/\alpha$.
If $|\alpha|=2^a3^b$, then the iterated application of the previous procedure to suitable powers of $\alpha$, produces a suitable $|\alpha|$-divisible set of curves on a surface $Y$, which is a smooth surface birational to $X/\alpha$.

\subsection{K3 surfaces}\label{subsec: preliminaries on K3}
We work with smooth projective complex surface.
\begin{definition}
A surface $Y$ is called K3 surface if its canonical bundle is trivial and $h^{(1,0)}(Y)=0$.
\end{definition}

The second cohomology group of a K3 surface equipped with the cup product is the unique even unimodular lattice of rank 22 and  signature $(3,19)$ and it is denoted by $\Lambda_{K3}$. The N\'eron--Severi group of a K3 surface $Y$ is a primitively embedded sublattice of $\Lambda_{K3}$ with signature $(1,\rho(Y)-1)$. Consequently, the transcendental lattice, which is the orthogonal to the N\'eron--Severi group in the second cohomology group, is a primitively embedded sublattice of $\Lambda_{K3}$ with signature $(2,20-\rho(Y))$.

Let $G\subset\Aut(Y)$ be a group of automorphisms of $Y$. We will say that it acts symplectically if it preserves the symplectic structure of $Y$, i.e. if its action on $H^{2,0}(Y)$ is trivial. 

The finite groups acting symplectically on a K3 surface are classified by Nikulin, \cite{Niksympl}, in the case of the Abelian group, and by Mukai, \cite{M}, in the other cases. A complete list can be found in \cite{X}.

If a finite group $G$ acts symplectically on a K3 surface $Y$, then $Y/G$ is a singular surface, whose desingularization $\widetilde{Y/G}$ is a K3 surface. 
\begin{definition} Let $Y$ be a K3 surface admitting a symplectic action of a finite group $G$. Let $\widetilde{Y/G}$ be the minimal model of $Y/G$. We will denote by $E_G$ the sublattice of $NS(\widetilde{Y/G})$ generated by the curves arising from the desingularization of $Y/G$. We will denote by $M_G$ the minimal primitive sublattice of $NS(\widetilde{Y/G})$ which contains $E_G$. We observe that $M_G$ is an overlattice of finite index of $E_G$. \end{definition}

We show now an explicit and very classic example: let $Y$ be a K3 surface which admits a symplectic action of $\Z/2\Z$. Then $Y/(\Z/2\Z)$ has 8 singular points of type $A_1$. Then the desingularization of $Y/G$ introduces 8 rational curves on $\widetilde{Y/(\Z/2\Z)}$, let us denote them by $M_i$, $i=1,\ldots, 8$. The lattice spanned by the curves $M_i$ is clearly isomorphic to $A_1^8$, so $E_{\Z/2\Z}=A_1^8$. 

One can also consider a different construction: one blows up $Y$ in the eight fixed points for the action of $\Z/2\Z$. One obtains the surface $\widetilde{Y}$, with 8 exceptional curves $E_i$, $i=1,\ldots, 8$. Then one lifts the action of $\Z/2\Z$ on $Y$ to an action of $\Z/2\Z$ on $\widetilde{Y}$, which fixes the exceptional curves. So one obtains a smooth surface $\widetilde{Y}/(\Z/2\Z)$, which is in fact isomorphic to $\widetilde{Y/(\Z/2\Z)}$. The $2:1$ map $\widetilde{Y}\ra \widetilde{Y}/(\Z/2\Z)$ is ramified on the union of the curves $E_i$ and so it is branched along the union of the curves $M_i$. By Section \ref{sec: covers}, it follows that $\sum_i M_i$ is divisible by 2 in $NS(\widetilde{Y/G})$. Hence $E_{\Z/2\Z}$ is generated by $M_i$, $i=1,\ldots 8$ and it is isometric to $A_1^8$; $M_{\Z/2\Z}$ is generated by the same classes as $E_{\Z/2\Z}$ and by the divisible class $\sum_i M_i/2$.

Similarly one can apply the results of section \ref{sec: covers} to the cyclic groups of order 3,4,6 and 8 in order to conclude that the K3 surface $\widetilde{Y/G}$ contains a divisible set of rational curves. The same is true also for cyclic groups of order 5 and 7 as proved by Nikulin in \cite{Niksympl}. This shows that for every cyclic group $G$ acting symplectic on a K3 surface $Y$, there is a $|G|$-divisible set of rational curves on the minimal model of $Y/G$. The description of this $|G|$-divisible set is given in \cite{Niksympl} and implies the description of the lattice $M_G$.

Let us assume that the sum of $n$ disjoint rational curves is divisible by 2 in $NS(Y)$. By Section \ref{sec: covers}, there exists a $2:1$ cover of $Y$ branched along the union of these curves. The covering surface is not minimal but one can contract certain curves in order to obtain a minimal model.  It is proved by Nikulin that only two possibilities occur: the minimal model of the covering surface is a K3 surface and in this case $n=8$ or the minimal model of the covering surface is an Abelian surface and in this case $n=16$. 
A similar results holds for covers of order 3 and is proved by Barth, \cite{Barthnine}. We collect these results in Proposition \ref{prop: classes divisible by 2 and 3} after introducing some definitions.

\begin{definition} An $A_k$ (resp. $D_m$, $m\geq 4$, $E_l$, $l=6,7,8$) configuration of curves is a set of $k$ (resp. $m$, $l$) irreducible smooth rational curves whose dual diagram is a Dynkin diagram of type $A_k$ (resp. $D_m$, $E_l$). 

A set of disjoint $A_k$ configurations is $n$ divisible if there is a linear combination of the curves contained in the configuration which can be divided by $n$ in the N\'eron--Severi group of the surface. 
\end{definition}

\begin{proposition}\label{prop: classes divisible by 2 and 3}{\rm(See \cite{NikKummer} and \cite{Barthnine})} Let $Y$ be a K3 surface which contains a set of $m$ disjoint rational curves (i.e. a set of $m$ disjoint $A_1$-configurations). If this set is divisible by 2, then $m$ is either 8 or 16. In the first case the cover surface associated to the divisible class is a K3 surface, in the latter it is an Abelian surface.

Let $Y$ be a K3 surface which contains a set of $m$ disjoint $A_2$-configurations. If this set is divisible by 3, then $m$ is either 6 or 9. In the first case the cover surface associated to the divisible class is a K3 surface, in the latter it is an Abelian surface. \end{proposition}

In case $G$ is an Abelian group acting symplectically on a K3 surface $Y$, the type and the number of points with a non trivial stabilizer is determined by Nikulin in \cite[Section 5]{Niksympl}. In the same paper the author determines the lattice $E_G$ and $M_G$ for all the admissible Abelian groups (we will recall this result in Proposition \ref{prop: EG and MG}). 

In certain cases the presence of certain configurations of rational curves suffices to conclude that the K3 surface is covered either by an Abelian or by a K3 surface. Since this property will be useful, we summarize the cases where it appears.
\begin{proposition}\label{prop: known results on covers of K3 with curves}
Let $Y$ be a K3 surface which admits 16 disjoint rational curves. Then it is the desingularization of the quotient of an Abelian surface by the group $\Z/2\Z$. In particular the set of these 16 curves is 2 divisible {\rm(see \cite[Theorem 1]{NikKummer} }.

Let $Y$ be a K3 surface which admits 15 disjoint rational curves. Then it is the desingularization of the quotient of a K3 surface by the group $(\Z/2\Z)^{4}$. In particular the set of these 15 rational curves contains 4 independent subsets of 8 rational curves which are 2-divisible {\rm (see \cite[Theorem]{Barthnine})}.

Let $Y$ be a K3 surface which admits 14 disjoint rational curves. Then it is the desingularization of the quotient of a K3 surface by the group $(\Z/2\Z)^{3}$. In particular the set of these 14 rational curves contains 3 independent subsets of 8 rational curves which are 2-divisible {\rm (see \cite[Corollary 8.7]{GS})}.

Let $Y$ be a K3 surface which admits 9 disjoint $A_2$-configurations of rational curves. Then it is the desingularization of the quotient of an Abelian surface by the group $\Z/3\Z$. In particular the set of 9 disjoint $A_2$-configurations of rational curves is 3 divisible  {\rm (see \cite[Theorem 8.6]{GS})}. 
\end{proposition}

\section{A preliminary and fundamental result}\label{sec: Bertin lemma and implications}
In this section we recall a result by Bertin, \cite[Lemma 3.1]{B}
and we deduce a corollary of this result, Proposition \ref{prop: the roots of L are smooth  irreducible curves}. These results are essential for the following. 

First we introduce some notation,
following \cite{B}. Let us consider a K3 surface $Y$. We denote by $\mathcal{C}^+$ the component of the cone $\{v\in
NS(Y)\otimes \R\mbox{ such that } v^2>0\}$ which contains at least
one ample class. We observe that all the ample cone is
contained in $\mathcal{C}^+$. Let us denote by
$$\Delta:=\{\delta\in NS(Y)\mbox{ such that }\delta^2=-2\}\mbox{ and by }\Delta^+:=\{\delta\in \Delta\mbox{ such that }\delta\mbox{ is effective}\}.$$
Moreover, we denote by $$B:=\{b\in \Delta^+\mbox{ such that
}b\mbox{ is the class of an irreducible curve}\}.$$ We observe
that the curves $C$, whose classes are contained in the
set $B$ are smooth irreducible rational curves. We pose
$$\mathcal{K}:=\{v\in \mathcal{C}^+\mbox{ such that }vb>0\mbox{ for all }b\in B\}.$$
The cone $\mathcal{K}$ is the ample cone of $S$ and
so $NS(K)\cap\overline{\mathcal{K}}\cap \mathcal{C}^+$ is the set of the
pseudo ample divisors of $Y$. This means that if $h\in
NS(Y)\cap\overline{\mathcal{K}}\cap \mathcal{C}^+$, then $h^2>0$ and $hv\geq 0$ for all the effective classes $v$.
\begin{lemma}{\rm (\cite[Lemma 3.1]{B})}\label{lemma: by Bertin} Let $Y$ be a projective K3 surface and $h\in
NS(Y)\cap\overline{\mathcal{K}}\cap \mathcal{C}^+$. Let us denote
$\Delta_h:=\Delta\cap h^{\perp}$. Then $B\cap h^{\perp}$ is a
basis of $\Delta_h$.\end{lemma} 

\begin{proposition}\label{prop: the roots of L are smooth irreducible curves}
Let $Y$ be a K3 surface. Let $h$ be a pseudoample divisor on $Y$
and let $L=h^{\perp}:=\{l\in NS(Y)\mbox{ \rm such that } lh=0\}$ be the orthogonal of $h$ in $NS(Y)$. Let us assume that there exists  a root lattice $R$ such that:\\
$(1)$ $L$ is an overlattice of finite index of $R$\\
$(2)$ the roots of $R$ and of $L$ coincide.\\ Then there exists a
basis of $R$ which is supported on smooth irreducible rational
curves.
\end{proposition}
\proof Let us consider the root lattice of $L$, denoted by $R(L)$. By definition of
$\Delta_h$, $R(L)$ and $\Delta_h$ coincide. So there is a basis
of $R(L)$ which is supported on smooth irreducible rational
curves, by \cite[Lemma 3.1]{B} (i.e. by Lemma \ref{lemma: by Bertin}).
By the hypothesis $R(L)$ is isometric to $R$ and so there exists a
basis for $R$ which is supported on smooth irreducible rational
curves.\endproof

\section{K3 surfaces quotients of Abelian surfaces}\label{sec: K3 quotients of Abelian surfaces}
In this section we concentrate on K3 surfaces which are constructed as quotients of an Abelian surface by a group of finite order. First we recall some known results and we compute the lattices associated to this construction in case $G$ is not Abelian. The results about these lattices are summarized in Proposition \ref{prop: definition lattices FG and KG}. Then we state and prove the main results of this section, which are the theorems \ref{theorem: X is A/G iff KG is in NS(X)} and \ref{theor: it suffices to have the curves on X}.
 
\subsection{Preliminaries and known results}

\begin{definition}
Let $A$ be an Abelian surface. Let $G\subset Aut(A)$ be a finite
group of automorphisms of $A$. Let us consider the minimal model
of $A/G$ and let us call it $X_G$.

Let $K_i$ be the curves on $X_G$ arising by the resolution of the
singularities of $A/G$. Let $F_G$ be the lattice spanned by the
curves $K_i$ and let $K_G$ be the minimal primitive sublattice of
$NS(X_G)$ containing the curves $K_i$. Clearly $K_G$ is an
overlattice of finite index, $r_G$, of $F_G$. We will say that the lattice $K_G$ is a lattice of Kummer type.
\end{definition}

The following well known result, due to Fujiki, classifies the
group $G\subset \Aut(A)$ such that $G$ does not contain
translations and $X_G$ is a K3 surface:

\begin{theorem}{\rm (\cite{F})} Let $G$ be a group of
automorphisms of an Abelian surface $A$ which does not contain
translations. If the minimal resolution of $A/G$ is a K3 surface,
then $G=\Z/n\Z$, $n=2,3,4,6$ or $G\in\{\mathbb{D}_8,
\mathbb{D}_{12}, \mathbb{T}\}$.
\end{theorem}

The requirement that $G$ does not contain translations is not  seriously restrictive, indeed the quotient of an Abelian surface by a finite group of translations produces another Abelian surface. Up to replacing the first Abelian surface by its quotient by translations, we can assume without loss of generalities that the group $G$ does not contain translations.

\subsection{Non-cyclic quotients of Abelian surfaces}\label{subsection: non cyclic quotients of A}
The aim of this Section is to describe the lattices $K_G$ in case $G$ is not cyclic. This lattices were computed also in \cite{W}, but two of the results given in \cite[Proposition 2.1]{W} are wrong. In particular we prove that the lattices $K_G$ are not the ones given in \cite[Proposition 2.1]{W} if $G=\mathbb{D}_8'$ and $G=\mathbb{D}_{12}$. 
\subsubsection{The actions of $\mathbb{D}_8$ and $\mathbb{D}_8'$} Let $G$ be the quaternion group. There are two different families of tori on which we can define the action of $G$ in such a way that $X_G$ is a K3 surface and on these two different families $G$ has
different sets of points with non trivial stabilizer, so the quotients of an Abelian surface by each of these actions 
produce two different singular surfaces with different
sets and types of singularities. 

The group $G$ has the
following presentation
$$\langle \alpha_4, \beta|\alpha_4^4=\beta^4=1\ \ \alpha_4^2=\beta^2,\ \alpha_4^{-1}\beta\alpha_4=\beta^{-1}\rangle.$$

We pose $A:=\C^2/\Lambda$ and $$\begin{array}{ll}\alpha_4:A\ra A,\ \ &(z_1,z_2)\mapsto (iz_1, -iz_2)\\

\beta:A\ra A,\ \ &(z_1,z_2)\mapsto (-z_2,
z_1)\end{array}
$$

The action of $\alpha_4$ and $\beta$ is algebraic as automorphism both of $A:=E_i\times E_i$ (where $E_i$ is the elliptic curve with
$j$-invariant equal to 1728, i.e. it is the elliptic curve
associated to the lattice $\langle 1,i\rangle$) and on the Abelian
surface $A':=\C^2/\Lambda$, where $\Lambda:=\langle (1,0), (i,0),
(\frac{1+i}{2},\frac{1+i}{2}),
(\frac{1+i}{2},\frac{i-1}{2})\rangle$. So the group generated by
$\alpha_4$ and $\beta$ is both a subgroup of $\Aut(A)$ and of
$\Aut(A')$. We denote it as $\mathbb{D}_8$ when it is considered
as subgroup of $\Aut(A)$ and as $\mathbb{D}_8'$ when it is
considered as subgroup of $\Aut(A')$.

We now identify the points of $A$ (resp. $A'$) which have a non
trivial stabilizer for $\mathbb{D}_8$ (resp. $\mathbb{D}_8$). All
of them are 2-torsion points and indeed are fixed for
$\alpha_4^2=\beta^2$. We have the following tables:
\begin{align}\label{eq: stable points D8}\begin{array}{c}
\mbox{ Points of $A$ with non trivial stabilizer for
$\mathbb{D}_8$}\\
\hline
\begin{array}{|c|ccc|}points\ in\ the\ same\ orbit & stabilizer&&\\
\hline
(0,0)&\mathbb{D}_8&=&\langle \alpha_4,\beta\rangle\\
\hline (\frac{1+i}{2},\frac{1+i}{2})&\mathbb{D}_8&=&\langle
\alpha_4,\beta\rangle\\ \hline (\frac{1}{2},\frac{1}{2}),
(\frac{i}{2},\frac{i}{2})&\Z/4\Z&=&\langle\beta\rangle\\ \hline (0,\frac{1+i}{2}),
(\frac{1+i}{2},0)&\Z/4\Z&=&\langle\alpha_4\rangle\\ \hline
(\frac{1}{2},\frac{i}{2}),
(\frac{i}{2},\frac{1}{2})&\Z/4\Z&=&\langle\alpha_4\circ\beta\rangle\\
\hline(\frac{1}{2},0),\ (\frac{i}{2}, 0),\ (0,\frac{1}{2}),\
(0,\frac{i}{2})&\Z/2\Z&=&\langle\alpha_4^2\rangle=\langle\beta^2\rangle\\
\hline(\frac{1+i}{2},\frac{1}{2}),\ (\frac{1+i}{2},\frac{i}{2}),\
(\frac{1}{2},\frac{1+i}{2}),\
(\frac{i}{2},\frac{1+i}{2})&\Z/2\Z&=&\langle\alpha_4^2\rangle=\langle\beta^2\rangle\\\hline\end{array}\end{array}
\end{align}

\begin{eqnarray}\label{table: non trivial stab for D8'}\end{eqnarray}
\vspace{-1cm}
$$
\begin{array}{c}
\mbox{ Points of $A'$ with non trivial stabilizer for
$\mathbb{D}_8'$ }\\
\hline
\begin{array}{|c|ccc|}points\ in\ the\ same\ orbit & stabilizer&&\\
\hline
(0,0)&\mathbb{D}_8'&=&\langle \alpha_4,\beta\rangle\\
\hline (\frac{1+i}{2},0)&\mathbb{D}_8'&=&\langle
\alpha_4,\beta\rangle\\
\hline (\frac{i}{2},\frac{i}{2})&\mathbb{D}_8'&=&\langle\alpha_4,\beta\rangle\\
\hline (\frac{1}{2},\frac{i}{2})&\mathbb{D}_8'&=&\langle\alpha_4,\beta\rangle\\
\hline (\frac{1}{2},0),(\frac{i}{2},0),(0,\frac{i}{2}),
(\frac{1+i}{2},\frac{i}{2})&\Z/2\Z&=&\langle\alpha_4^2\rangle=\langle\beta^2\rangle\\
\hline
(\frac{1+i}{4},\frac{1+i}{4}),\ (\frac{1-i}{4},\frac{i-1}{4}),\
(\frac{i-1}{4},\frac{i-1}{4}), \
(\frac{-1-i}{4},\frac{i+1}{4})&\Z/2\Z&=&\langle\alpha_4^2\rangle=\langle\beta^2\rangle\\\hline

(\frac{1+i}{4},\frac{i-1}{4}),\ (\frac{1-i}{4},\frac{i+1}{4}),\
(\frac{i-1}{4},\frac{i+1}{4}), \
(\frac{-1-i}{4},\frac{i-1}{4})&\Z/2\Z&=&\langle\alpha_4^2\rangle=\langle\beta^2\rangle\\\hline\end{array}\end{array}
$$

We observe that $\alpha_4^2(=\beta^2 )$ is the center of
$G:=\langle\alpha_4,\beta\rangle$ and in particular it is a normal
subgroup of $\langle\alpha_4,\beta\rangle$. So, in order to
construct $A/G$, one can first consider $A/\langle
\alpha_4^2\rangle$ and then one can consider $(A/\langle
\alpha_4^2\rangle)/Q$ where $Q$ is the quotient group
$G/\langle\alpha_4^2\rangle$. The group $Q$ is isomorphic to
$(\Z/2\Z)^2$ and is generated by $\overline{\alpha_4}$ and
$\overline{\beta}$ where $\overline{g}$ is the image of $g\in G$
under the quotient map $G\ra G/\langle\alpha_4^2\rangle$.

The automorphism $\alpha_4^2$ is $(z_1,z_2)\ra (-z_1,-z_2)$ and the
surface $A/\alpha_4^2$ is singular in the image of the 16
2-torsion points of $A$. We denote by $A[2]$ the set of these points on $A$. The desingularization, $Km(A)$, of $A/\langle \alpha_4^2\rangle$ is obtained by blowing up once the
singular points and it is the Kummer surface of $A$. Let $(p,q)\in
A[2]$, then we denote by $K_{(p,q)}\subset Km(A)$ the 16 rational
curves arising from the desingularization of
$A/\langle\alpha_4^2\rangle$.

The minimal resolution $X_G$ of $A/G$ is birational
to the minimal resolution of $Km(A)/Q$ where
$Q=G/\langle \alpha_4^2\rangle$. Since the minimal
resolution of $A/G$ is a K3 surface, and birational K3
surfaces are isomorphic, we conclude that the minimal model
of $A/G$ is the minimal model of $Km(A)/Q$. We recall that
$Km(A)$ is a K3 surface and $Q=(\Z/2\Z)^2$ acts on $Km(A)$
preserving the symplectic structure (indeed the quotient has the
induced symplectic structure). The action of the group
$(\Z/2\Z)^2$ on a K3 surface is very well known, see \cite[Section 5]{Niksympl}: each copy of $\Z/2\Z$ in $(\Z/2\Z)^2$ stabilizes exactly 8 points and there are no points fixed by the full group. So there are 24 points with non trivial stabilizer in $Km(A)$ and then in the quotient $Km(A)/Q$ we have 12 singular points, each of them of type $A_1$. This can be also checked by hand, considering the action of $Q$ over the curves $K_{(p,q)}$. By Tables \eqref{eq: stable points D8}, $G$ has no points with a non trivial stabilizer outside the set of the 2-torsion points $A[2]$. So the points with a non trivial stabilizer for $Q$ on $Km(A)$ are all contained in the curves $K_{(p,q)}$.

If the point $(p,q)\in A$ is fixed by $G$, then
$Q\simeq (\Z/2\Z)^2$ preserves the curve $K_{(p,q)}$, which is a
copy of $\mathbb{P}^1$. Hence $Q$ has six points with
non trivial stabilizer on $K_{(p,q)}$. So in $Km(A)/Q$ there are
three singular points on the image of $K_{(p,q)}$. Hence on
$X_G$ there are three rational curves (arising from the
desingularization of these points) which intersect the image of
$K_{(p,q)}$. We call these curves $K_{(p,q)}^{(i)}$, $i=1,2,3$ and
$K_{(p,q)}^{(0)}$ the strict transform of the image of $K_{(p,q)}$ under the quotient map
$Km(A)\ra Km(A)/Q$. The curves $K_{(p,q)}^{(i)}$, $i=0,1,2,3$
generate a copy of the lattice $D_4$ in $NS(X_G)$.

If the point $(p,q)\in A$ has the group $\Z/4\Z\subset G$ as stabilizer,
then there is another point $(p',q')$ in its orbit. The quotient
group $Q$ switches the curves $K_{(p,q)}$ and $K_{(p',q')}$ and
has 4 points with non trivial stabilizer in $K_{(p,q)}\cup
K_{(p',q')}$ (2 on each curve). So on the quotient surface
$Km(A)/Q$ there is a curve $K_{(p,q)}^{(0)}$ which is the common
image of $K_{(p,q)}$ and of $K_{(p',q')}$ and there are 2 singular
points on such a curve. We denote by $K_{(p,q)}^{(1)}$ and
$K_{(p,q)}^{(2)}$ the curves in $X_G$ arising from the
desingularization of these two singular points. The curves
$K_{(p,q)}^{(i)}$, $i=0,1,2$ generate a copy of $A_3$ in $NS(X_G)$ (here, with an abuse of notation, we denote by $K_{(p,q)}^{(0)}$ both a curve on $Km(A)/Q$ and its strict transform on $X_G$).

If the point $(p,q)\in A$ has the group $\Z/2\Z\subset G$ as stabilizer, it
is generated by $\alpha_4^2$ and in the same orbit of $(p,q)$
there are the other 3 points $(p',q')$, $(p'',q'')$,
$(p''',q''')$. The group $Q$ permutes the curves $K_{(p,q)}$,
$K_{(p',q')}$, $K_{(p'',q'')}$, $K_{(p''',q''')}$ in $Km(A)$. So
their image in $X_G$ is a unique curve $K_{(p,q)}$.

The curves arising from the desingularization of $A/\mathbb{D}_8$
(resp. $A'/\mathbb{D}_8'$) are represented in Figure \ref{figure:
A/D8 curves} (resp. Figure \ref{figure: A'/D8' curves}).

\begin{figure}\scalebox{0.45}{
\begin{tikzpicture}
\node at (-0.5 ,0) {$K_{(0,0)}^{(0)}$};

\node at (1,1.5) {$K_{(0,0)}^{(1)}$};

\node at (2,-1) {$K_{(0,0)}^{(2)}$};

\node at (3,1.5) {$K_{(0,0)}^{(3)}$};
    \draw (0,0) -- (4,0);
    \draw (1,1)--(1,-0.5);
    \draw (2,1)--(2,-0.5);
    \draw (3,1)--(3,-0.5);

\node at (5.5 ,0) {$K_{(\frac{1+i}{2},\frac{1+i}{2})}^{(0)}$};

\node at (7,1.5) {$K_{(\frac{1+i}{2},\frac{1+i}{2})}^{(1)}$};

\node at (8,-1) {$K_{(\frac{1+i}{2},\frac{1+i}{2})}^{(2)}$};

\node at (9,1.5) {$K_{(\frac{1+i}{2},\frac{1+i}{2})}^{(3)}$};
    \draw (6,0) -- (10,0);
    \draw (7,1)--(7,-0.5);
    \draw (8,1)--(8,-0.5);
    \draw (9,1)--(9,-0.5);

\node at (11.5 ,0) {$K_{(\frac{1}{2},\frac{1}{2})}^{(0)}$};

\node at (13,1.5) {$K_{(\frac{1}{2},\frac{1}{2})}^{(1)}$};

\node at (14,-1) {$K_{(\frac{1}{2},\frac{1}{2})}^{(2)}$};
    \draw (12,0) -- (15,0);
    \draw (13,1)--(13,-0.5);
    \draw (14,1)--(14,-0.5);

\node at (16.5 ,0) {$K_{(0,\frac{1+i}{2})}^{(0)}$};

\node at (18,1.5) {$K_{(0,\frac{1+i}{2})}^{(1)}$};

\node at (19,-1) {$K_{(0,\frac{1+i}{2})}^{(2)}$};
    \draw (17,0) -- (20,0);
    \draw (18,1)--(18,-0.5);
    \draw (19,1)--(19,-0.5);

\node at (21.5 ,0) {$K_{(\frac{1}{2},\frac{i}{2})}^{(0)}$};

\node at (23,1.5) {$K_{(\frac{1}{2},\frac{i}{2})}^{(1)}$};

\node at (24,-1) {$K_{(\frac{1}{2},\frac{i}{2})}^{(2)}$};
    \draw (22,0) -- (25,0);
    \draw (23,1)--(23,-0.5);
    \draw (24,1)--(24,-0.5);

\node at (26 ,1.5) {$K_{(\frac{1}{2},0)}^{(0)}$};

\draw (26.5,1)--(26.5,-1);

\node at (27.5,1.5) {$K_{(\frac{1+i}{2},\frac{1}{2})}^{(0)}$};
    \draw (28,1)--(28,-1);

\end{tikzpicture}
} \caption{Curves of $F_{\mathbb{D}_8}$ on
$X_{\mathbb{D}_8}$}\label{figure: A/D8 curves}
\end{figure}

\begin{figure}\scalebox{0.45}{
\begin{tikzpicture}
\node at (-0.5 ,0) {$K_{(0,0)}^{(0)}$};

\node at (1,1.5) {$K_{(0,0)}^{(1)}$};

\node at (2,-1) {$K_{(0,0)}^{(2)}$};

\node at (3,1.5) {$K_{(0,0)}^{(3)}$};
    \draw (0,0) -- (4,0);
    \draw (1,1)--(1,-0.5);
    \draw (2,1)--(2,-0.5);
    \draw (3,1)--(3,-0.5);

\node at (5.5 ,0) {$K_{(\frac{1+i}{2},0)}^{(0)}$};

\node at (7,1.5) {$K_{(\frac{1+i}{2},0)}^{(1)}$};

\node at (8,-1) {$K_{(\frac{1+i}{2},0)}^{(2)}$};

\node at (9,1.5) {$K_{(\frac{1+i}{2},0)}^{(3)}$};
    \draw (6,0) -- (10,0);
    \draw (7,1)--(7,-0.5);
    \draw (8,1)--(8,-0.5);
    \draw (9,1)--(9,-0.5);

\node at (11.5 ,0) {$K_{(\frac{i}{2},\frac{i}{2})}^{(0)}$};

\node at (13,1.5) {$K_{(\frac{i}{2},\frac{i}{2})}^{(1)}$};

\node at (14,-1) {$K_{(\frac{i}{2},\frac{i}{2})}^{(2)}$};

\node at (15,1.5) {$K_{(\frac{i}{2},\frac{i}{2})}^{(3)}$};

    \draw (12,0) -- (16,0);
    \draw (13,1)--(13,-0.5);
    \draw (14,1)--(14,-0.5);
\draw (15,1)--(15,-0.5);

\node at (17.5 ,0) {$K_{(\frac{1}{2},\frac{i}{2})}^{(0)}$};

\node at (19,1.5) {$K_{(\frac{1}{2},\frac{i}{2})}^{(1)}$};

\node at (20,-1) {$K_{(\frac{1}{2},\frac{i}{2})}^{(2)}$};

\node at (21,1,5) {$K_{(\frac{1}{2},\frac{i}{2})}^{(3)}$};

    \draw (18,0) -- (22,0);

    \draw (19,1)--(19,-0.5);
    \draw (20,1)--(20,-0.5);
    \draw (21,1)--(21,-0.5);

\node at (23 ,1.5) {$K_{(\frac{1}{2},0)}^{(0)}$};

\draw (23,1)--(23,-1);

\node at (24.5,1.5) {$K_{(\frac{1+i}{4},\frac{1+i}{4})}^{(0)}$};
    \draw (25,1)--(25,-1);

\node at (26,1.5) {$K_{(\frac{1+i}{4},\frac{i-1}{4})}^{(0)}$};
    \draw (26.5,1)--(26.5,-1);

\end{tikzpicture}
}\caption{Curves of $F_{\mathbb{D}_8'}$ on
$X_{\mathbb{D}_8}'$}\label{figure: A'/D8' curves}
\end{figure}

\subsubsection{The lattice $K_{\mathbb{D}_8}$}
Let us now fix a specific action of $G$. In particular let the Abelian surface be $A$ (e.g. $A\simeq E_i\times E_i$) and so $G\subset \Aut(A)$ is the group $\mathbb{D}_8$.
In this case the lattice $F_{\mathbb{D}_8}$ is isometric to
$D_4^2\oplus A_3^3\oplus A_1^2$, see Table \eqref{table: non trivial stab for  D8'}. Its discriminant group is
$(\Z/2\Z)^6\times (\Z/4\Z)^3$ and is generated by the following
classes:
$d_1:=\frac{1}{2}\left(K_{(0,0)}^{(1)}+K_{(0,0)}^{(2)}\right)$,
$d_2:=\frac{1}{2}\left(K_{(0,0)}^{(1)}+K_{(0,0)}^{(3)}\right)$,
$d_3:=\frac{1}{2}\left(K_{(\frac{1+i}{2},\frac{1+i}{2})}^{(1)}+K_{(\frac{1+i}{2},\frac{1+i}{2})}^{(2)}\right)$,
$d_4:=\frac{1}{2}\left(K_{(\frac{1+i}{2},\frac{1+i}{2})}^{(1)}+K_{(\frac{1+i}{2},\frac{1+i}{2})}^{(3)}\right)$,
$d_5:=\frac{1}{2}K_{(\frac{1}{2},0)}$,
$d_6:=\frac{1}{2}K_{(\frac{1+i}{2},\frac{1}{2})}$,\\
$d_7:=\frac{1}{4}\left(K_{(\frac{1}{2},\frac{1}{2})}^{(1)}+2K_{(\frac{1}{2},\frac{1}{2})}^{(0)}+3K_{(\frac{1}{2},\frac{1}{2})}^{(2)}\right)$,
$d_8:=\frac{1}{4}\left(K_{(0,\frac{1+i}{2})}^{(1)}+2K_{(0,\frac{1+i}{2})}^{(0)}+3K_{(0,\frac{1+i}{2})}^{(2)}\right)$,
$d_9:=\frac{1}{4}\left(K_{(\frac{1}{2},\frac{i}{2})}^{(1)}+2K_{(\frac{1}{2},\frac{i}{2})}^{(0)}+3K_{(\frac{1}{2},\frac{i}{2})}^{(2)}\right)$.

The set of 12 curves $\mathcal{S}:=\{K_{(0,0)}^{(j)}$ , $K_{(\frac{1+i}{2},\frac{1+i}{2}))}^{(j)}$, $K_{(\frac{1}{2},\frac{1}{2})}^{(k)}$, $K_{(0,\frac{1+i}{2})}^{(k)}$, $K_{(\frac{1}{2},\frac{i}{2})}^{(k)}\}$, $j=1,2,3$, $k=1,2$ arises from the desingularization of the quotient of a K3 surface (the surface $Km(A)$) by the group $(\Z/2\Z)^2$. By \cite[Section 6, case 2a), equation (6.17)]{Niksympl}, $\mathcal{S}$ contains 2 independent subsets of 8 curves which are two divisible. Indeed the two classes $$v_1:=d_1+d_3+2d_7+2d_8,\ \ v_2:=d_2+d_4+2d_8+2d_9$$ 
are contained in $NS(X_G)$. Moreover, the set $\mathcal{S}\cup\{ K_{(\frac{1+i}{2},0)}, K_{(\frac{1+i}{2},\frac{1}{2})}\}$ forms a set of 14 disjoint rational curves contained in the
curves of $F_{\mathbb{D}_8}$ (this set consists of the vertical
curves in Figure \ref{figure: A/D8 curves}). By Proposition \ref{prop: known results on covers of K3 with  curves}, the minimal primitive sublattice of the N\'eron--Severi group which contains these 14 curves is spanned by the curves and by 3 other divisible classes. So there is another divisible class contained in $NS(X_G)$, which is: $$v_3:=d_1+d_3+2d_9+d_5+d_6.$$

Let us denote by $L_{\mathbb{D}_8}$ the lattice spanned by
$F_{\mathbb{D}_8}$ and by the classes $v_1$, $v_2$, $v_3$. Its
discriminant group is
$(\Z/4\Z)^3$ and is generated by:\\
$\delta_1:=d_4+d_7=\frac{1}{4}\left(2K_{(\frac{1+i}{2},\frac{1+i}{2})}^{(1)}+2K_{(\frac{1+i}{2},\frac{1+i}{2})}^{(2)}+K_{(\frac{1}{2},\frac{1}{2})}^{(1)}+2K_{(\frac{1}{2},\frac{1}{2})}^{(0)}+3K_{(\frac{1}{2},\frac{1}{2})}^{(2)}\right)$,\\
$\delta_2:=d_3+d_4+d_8=\frac{1}{4}\left(2K_{(\frac{1+i}{2},\frac{1+i}{2})}^{(2)}+2K_{(\frac{1+i}{2},\frac{1+i}{2})}^{(3)}+K_{(0,\frac{1+i}{2})}^{(1)}+2K_{(0,\frac{1+i}{2})}^{(0)}+3K_{(0,\frac{1+i}{2})}^{(2)}\right)$,\\
$\delta_3:=d_3+d_5+d_9=\frac{1}{4}\left(
2K_{(\frac{1+i}{2},\frac{1+i}{2})}^{(1)}+2K_{(\frac{1+i}{2},\frac{1+i}{2})}^{(2)}+2K_{(\frac{1}{2},0)}+K_{(\frac{1}{2},\frac{i}{2})}^{(1)}+2K_{(\frac{1}{2},\frac{i}{2})}^{(0)}+3K_{(\frac{1}{2},\frac{i}{2})}^{(2)}\right)$

There are two possibilities, either $K_{\mathbb{D}_8}\simeq
L_{\mathbb{D}_8}$ or $K_{\mathbb{D}_8}$ is an overlattice of
finite index of $L_{\mathbb{D}_8}$, see Section \ref{subsec: lattice}. In the latter case $K_{\mathbb{D}_8}$ contains an element $w$ which is non trivial in the discriminant group of $L_{\mathbb{D}_8}$. So $w=\sum_{i=1}^3\alpha_i\delta_i$,
$\alpha_i\in \Z$ and $(\alpha_1,\alpha_2,\alpha_3)\not\equiv
(0,0,0)\mod 4$. If
$(\alpha_1,\alpha_2,\alpha_3)\equiv(2,2,2)\mod 4$, let $z:=w$,
otherwise let $z:=2w$. The element $z\in K_{\mathbb{D}_8}$
consists of the sum of certain disjoint rational curves divided by
2. These curves are chosen in\\
$\{K_{(\frac{1}{2},\frac{1}{2})}^{(1)},K_{(\frac{1}{2},\frac{1}{2})}^{(2)},K_{(0,\frac{1+i}{2})}^{(1)},K_{(0,\frac{1+i}{2})}^{(2)},K_{(\frac{1}{2},\frac{i}{2})}^{(1)},K_{(\frac{1}{2},\frac{i}{2})}^{(2)}\}$
(which are the ones which appear in $\delta_i$ with an odd
coefficient). By Proposition \ref{prop: classes divisible by 2 and  3} a set of at most 6 disjoint rational curves can not be divisible by 2, so $w$ can not exist. We conclude that
$L_{\mathbb{D}_8}\simeq K_{\mathbb{D}_8}$ and it is generated by
the generators of $F_{\mathbb{D}_8}$ and by $\{v_1,v_2,v_3\}$.

This result agrees with the one given in \cite[Proposition 2.1]{W}.\\

\subsubsection{The lattice $K_{\mathbb{D}_8'}$}
We now consider the Abelian surface $A'$, so $G\subset\Aut(A')$ is $\mathbb{D}_8'$. In this
case $F_{\mathbb{D}_8'}$ is $D_4^4\oplus A_1^3$. Its discriminant
group is $(\Z/2\Z)^{11}$ and is generated by
$d_1':=\frac{1}{2}\left(K_{(0,0)}^{(1)}+K_{(0,0)}^{(2)}\right)$,
$d_2':=\frac{1}{2}\left(K_{(0,0)}^{(1)}+K_{(0,0)}^{(3)}\right)$,
$d_3':=\frac{1}{2}\left(K_{(\frac{1+i}{2},0)}^{(1)}+K_{(\frac{1+i}{2},0)}^{(2)}\right)$,
$d_4':=\frac{1}{2}\left(K_{(\frac{1+i}{2},0)}^{(1)}+K_{(\frac{1+i}{2},0)}^{(3)}\right)$,
$d_5':=\frac{1}{2}\left(K_{(\frac{i}{2},\frac{i}{2})}^{(1)}+K_{(\frac{i}{2},\frac{i}{2})}^{(2)}\right)$,
$d_6':=\frac{1}{2}\left(K_{(\frac{i}{2},\frac{i}{2})}^{(1)}+K_{(\frac{i}{2},\frac{i}{2})}^{(3)}\right)$,
$d_7':=\frac{1}{2}\left(K_{(\frac{1}{2},\frac{i}{2})}^{(1)}+K_{(\frac{1}{2},\frac{i}{2})}^{(2)}\right)$,
$d_8':=\frac{1}{2}\left(K_{(\frac{1}{2},\frac{i}{2})}^{(1)}+K_{(\frac{1}{2},\frac{i}{2})}^{(3)}\right)$,
$d_9':=\frac{1}{2}K_{(\frac{1}{2},0)}$,
$d_{10}':=\frac{1}{2}K_{(\frac{1+i}{4},\frac{1+i}{4})}$,
$d_{11}':=\frac{1}{2}K_{(\frac{1+i}{2},\frac{i-1}{2})}^{(1)}$.

The set of 12 curves $\mathcal{S}:=\{K_{(0,0)}^{(j)}$ , $K_{(\frac{1+i}{2},0)}^{(j)}$, $K_{(\frac{i}{2},\frac{i}{2})}^{(j)}$, $K_{(\frac{1}{2},\frac{i}{2})}^{(j)}\}$, $j=1,2,3$, arises from the desingularization of the quotient of a K3 surface (the surface $Km(A')$) by the group $(\Z/2\Z)^2$ (i.e. the group $Q:=G/\langle \alpha_4^2\rangle$) so, by \cite[Section 6, case 2a), equation (6.17)]{Niksympl}, there are 2 divisible classes whose curves are in $\mathcal{S}$. Hence
\begin{eqnarray}\label{eq: v1' and v2'} v_1':=d_1'+d_3'+d_5'+d_7',\ v_2':=d_2'+d_4'+d_6'+d_8'\end{eqnarray} 
are contained in $NS(X_G)$. Moreover,
in the lattice $F_{\mathbb{D}_8'}$ it is easy to identify a set of
15 disjoint rational curves (the vertical ones in Figure \ref{figure:
A'/D8' curves}), which contains the set $\mathcal{S}$. By Proposition \ref{prop: known results on covers of K3 with  curves} the
minimal primitive sublattice of the N\'eron--Severi group which
contains these curves is spanned by the curves and by 4 other
divisible classes. Two of these divisible classes are $v_1'$ and $v_2'$, the others are:
\begin{eqnarray}\label{eq: v3' and v4'}\ \ \ v_3':=d_1'+d_3'+d_4'+d_6'+d_9'+d_{10}',\ v_4':=d_1'+d_4'+d_7'+d_8'+d_9'+d_{11}'.\end{eqnarray}
These 4 divisible classes are also contained in
$K_{\mathbb{D}_8'}$. 
Let us denote by $L_{\mathbb{D}_8'}$ the lattice spanned by
$F_{\mathbb{D}_8'}$ and by the classes $v_1'$, $v_2'$, $v_3'$,
$v_4'$. Its discriminant group is $(\Z/2\Z)^{3}$ and is generated
by:
$$\delta_1':=d_2'+d_3'+d_4'+d_5',\ \delta_2':=d_3'+d_4'+d_6'+d_7',\ \delta_3'':=d_4'+d_5'+d_6'+d_7'+d_{11}'.$$

If $K_{\mathbb{D}_8'}$ does not coincide with
$L_{\mathbb{D}_8'}$, then there is a vector $w$ which is non trivial in the discriminant group of $L_{\mathbb{D}_8'}$, and is not contained in $K_{\mathbb{D}_8'}$, by Section \ref{subsec: lattice}.
The curves which appear with a non trivial coefficient in
$\delta_1'$, $\delta_2'$ and $\delta_3'$ are all contained in the
set of 15 disjoint rational curves considered above. So if a
vector as $w$ exists it gives an overlattice of the lattice
spanned by 15 disjoint rational curves with index greater then
$2^4$ and contained in the N\'eron--Severi group of a K3 surface,
but this is impossible: indeed if we construct an overlattice of index $2^4$ of $A_1^{15}$, every 2-divisible set contains exactly 8 disjoint rational curves by Proposition \ref{prop: classes divisible by 2 and  3} and two divisible sets have exactly 4 curves in common. Let us denote by $e_i$ the 15 classes generating $A_1^{15}$. The first divisible set contains 8 classes, so up to permutation of the indices we can assume that it is $\mathcal{S}_1:=\{e_1,\ldots, e_8\}$. The second one contains 8 classes, four of them in common with $\mathcal{S}_1$, so we can assume that it is $\mathcal{S}_2:=\{e1_,\ldots, e_4, e_9,\ldots e_{12}\}$. Similarly the third can be chosen to be $\mathcal{S}_3:=\{e_1,e_2,e_5,e_6,e_9,e_{10},e_{13}, e_{14}\}$. This forces the fourth to be $\mathcal{S}_4:=\{e_1,e_3,e_5,e_7,e_9,e_{11}, e_{13}, e_{15}\}$. But now it is not possible to find another subset of $\{e_1,\ldots e_{15}\}$ which contains 8 elements and such that its intersection with each set $\mathcal{S}_i$ contains exactly 4 elements. 

We conclude that
$L_{\mathbb{D}_8'}$ coincides with $K_{\mathbb{D}_8'}$, which is
generated by the vectors in $F_{\mathbb{D}_8'}$ and by the four
vectors $v_1'$, $v_2'$, $v_3'$ and $v_4'$.

This result is different to the one given in \cite{W}. Indeed the lattice of Kummer type $\Pi_{\mathbb{D}_8'}$ described in \cite[Proposition 2.1]{W} contains a vector which consists of six disjoint rational curves divided by 2, which
is not possible by Proposition \ref{prop: classes divisible by 2 and  3}.

\subsubsection{The lattice $K_{\mathbb{T}}$}
Let us now consider the torus $A'$. There is an extra
automorphism, which is not contained in $\mathbb{D}_8'$ and which acts on $A'$, the automorphism $\gamma:(z_1,z_2)\ra
(\frac{i-1}{2}(z_1-z_2),\frac{-i-1}{2}(z_1+z_2))$. The automorphism
$\gamma$ has order 3 and the group
$\langle\alpha_4,\beta,\gamma\rangle$ is the binary tetrahedral group $\mathbb{T}$. It is the
semidirect product $\langle \gamma\rangle\ltimes \mathbb{D}_8'$.
In particular $\mathbb{D}_8'$ is a normal subgroup of $\mathbb{T}$
hence $A'/\mathbb{T}$ is birational to
$(A'/\mathbb{D}_8')/\langle\overline{\gamma}\rangle$, where
$\overline{\gamma}$ is the image of $\gamma$ under the quotient
map $\mathbb{T}\ra\mathbb{T}/\mathbb{D}_8'$. Hence the K3 surface
$X_{\mathbb{D}_8'}$, desingularization of $A'/\mathbb{D}_8'$,
admits a symplectic automorphism, $\gamma_X$, of order 3 induced
by $\overline{\gamma}$. The K3 surface $X_{\mathbb{T}}$,
desingularization of $A/\mathbb{T}$, is then isomorphic to the K3
which is the desingularization of
$X_{\mathbb{D}_8'}/\gamma_X$. In order to construct
$F_{\mathbb{T}}$, we consider the action of $\gamma_X$ on the
curves of $F_{\mathbb{D}_8'}$, see Figure \ref{figure: A'/D8' curves}: since
$\gamma((\frac{1+i}{2},0))=(\frac{1}{2},\frac{i}{2})$ and
$\gamma((\frac{1}{2},\frac{i}{2}))=(\frac{i}{2},\frac{i}{2})$, the
three copies of $D_4$, whose components are
$K_{(\frac{1+i}{2},0)}^{(j)}$,
$K_{(\frac{1}{2},\frac{i}{2})}^{(j)}$ and
$K_{(\frac{i}{2},\frac{i}{2})}^{(j)}$, $j=0,1,2,3$ are permuted by
$\gamma_X$. Hence these three copies of $D_4$ are identified on
$X_{\mathbb{T}}$ and correspond to a unique copy of $D_4$ on
$X_{\mathbb{T}}$. The same happens to the three copies of $A_1$,
which are permuted by $\gamma_X$ and thus give a unique copy of
$A_1$ on $X_{\mathbb{T}}$. Since $(0,0)$ is a fixed point for
$\gamma$, the automorphism $\gamma_X$ preserves the set of curves
$\{K_{(0,0)}^{(j)}\}$, $j=0,1,2,3$. Indeed $\gamma_X$ preserves
the curve $K_{(0,0)}^{(0)}$ and permutes the curves
$K_{(0,0)}^{(j)}$, $j=1,2,3$. So it is not the identity on
$K_{(0,0)}^{(0)}$ (since it moves the intersection points among
$K_{(0,0)}^{(0)}$ and $K_{(0,0)}^{(j)}$, $j=1,2,3$) and thus has
two fixed point on it. On the quotient these two points
corresponds to two singularities of type $A_2$. This gives 6
curves on $X_{\mathbb{T}}$ (one is the image of
$K_{(0,0)}^{(0)}$, one is the common image of $K_{(0,0)}^{(j)}$
for $j=1,2,3$, four come from the desingularization of the
two singular points of type $A_2$) and their dual graph is a copy
of $E_6$ (the image of $K_{(0,0)}^{(0)}$ intersects the image of
$K_{(0,0)}^{(j)}$ and one curve of each copy of the two $A_2$
arising from the desingularization).

We recall that a symplectic automorphism of order 3 on a K3
surface has exactly 6 fixed points. Since $\gamma_X$ fixes two
points on $K_{(0,0)}^{(0)}$ and has no fixed points on the
other curves of $F_{\mathbb{D}_8'}$, it necessarily fixes 4 points
in $X_{\mathbb{D}_8'}$ outsides the curves in
$F_{\mathbb{D}_8'}$, hence the desingularization $X_{\mathbb{T}}$
introduces 4 disjoint $A_2$-configurations. Thus, the lattice $F_{\mathbb{T}}$ is isometric to
$E_6\oplus D_4\oplus A_1\oplus A_2^4$. We fix the following
notation:
$$
\xymatrix{
e_3\ar@{-}[r]&e_2\ar@{-}[r]&e_0\ar@{-}[r]\ar@{-}[d]&e_4\ar@{-}[r]&e_5&&f_1\ar@{-}[r]&f_0\ar@{-}[r]\ar@{-}[d]&f_2\\
&&e_1&&&&&f_3}$$ where $\{e_j\}$ forms a basis of $E_6$ and
$\{f_j\}$ forms a basis of $D_4$. We denote by $a^{(1)}$ generator
of $A_1$ and by $a_j^{(h)}$, $j=1,2$, $h=1,2,3,4$, the basis of the $h$-th copy of $A_2$. A basis for the discriminant group of
$F_{\mathbb{T}}$ is given by
$d_1:=\frac{1}{3}(e_2+2e_3+e_4+2e_5)+\frac{1}{2}(f_1+f_2)$
$d_2:=\frac{1}{3}(a_1^{(2)}+2a_2^{(2)})+\frac{1}{2}(f_1+f_3)$
$d_3:=\frac{1}{3}(a_1^{(3)}+2a_2^{(3)})+\frac{1}{2}a^{(1)}$
$d_4:=\frac{1}{3}(a_1^{(4)}+2a_2^{(4)})$,
$d_5:=\frac{1}{3}(a_1^{(5)}+2a_2^{(5)})$.

The curves $e_2,e_3,e_4$, $e_5$, $a_1^{(j)}$, $a_2^{(j)}$, $j=2,3,4,5$ are the curves arising from the resolution of the quotient $X_{\mathbb{D}_8'}/\gamma_X$. So
by Section \ref{sec: covers} (see also Proposition \ref{prop: classes divisible by 2 and 3}), the class
$\left(e_2+2e_3+e_4+2e_5+\sum_{j=1}^4(a_1^j+2a_2^j)\right)/3$ is
contained in $NS(X_{\mathbb{T}})$ and hence also in
$K_{\mathbb{T}}$ (which is the minimal primitive sublattice of
$NS(X)_{\mathbb{T}}$ which contains the curves $e_h$, $f_j$,
$a_r^{(s)}$). So the vector $v:=4d_1+4d_2+4d_3+d_4+d_5\mod
F_{\mathbb{T}}$ is contained in $K_{\mathbb{T}}$. Let us denote by
$L_{\mathbb{T}}$ the lattice generated by the curves of
$F_{\mathbb{T}}$ and by $v$. Its discriminant group is generated
by $\delta_1:=d_1+d_2$, $\delta_2:=d_1+d_3$, $\delta_3:=d_1+d_4$.
If $L_{\mathbb{T}}\neq K_{\mathbb{T}}$, then there exists a vector
$w\in K_{\mathbb{T}}$ which is a non trivial element of the discriminant group of $L_{\mathbb{T}}$, which is $(\Z/6\Z)^3$. So either $w$ or a multiple of $w$
generates either $\Z/3\Z$ or $\Z/2\Z$ in the discriminant group of
$L_{\mathbb{T}}$. Every linear combination of $\delta_1$,
$\delta_2$ and $\delta_3$ which generates $\Z/2\Z$ is the sum of
at most 4 disjoint rational curves divided by 2 and so can not be a class in $NS(X_{\mathbb{T}})$, by Proposition \ref{prop: classes divisible by 2 and 3}.
Similarly, every linear combination of $\delta_1$, $\delta_2$ and
$\delta_3$ which generates $\Z/3\Z$ contains at most 5 disjoint $A_2$-configurations of rational curves. By Proposition \ref{prop: classes divisible by 2 and 3} it is impossible to construct a 3-divisible class with less then 6 disjoint $A_2$-configurations. We conclude that
$K_{\mathbb{T}}=L_{\mathbb{T}}$ is generated by $v$ and by the
curves in $F_{\mathbb{T}}$.

This result agrees with the one given in \cite{W}.

\subsubsection{The lattice $K_{\mathbb{D}_{12}}$}

Let $A$ be the Abelian surface $A:=E_{\zeta_3}\times E_{\zeta_3}$
where $\zeta_3$ is a primtive 3-rd root of unity and $E_{\zeta_3}$
is the elliptic curve with $j$-invariant 0. Let us now consider
the action of the group $\mathbb{D}_{12}$, which is algebraic on $A$
and is generated by the two automorphisms
$\alpha_6:(z_1,z_2)\mapsto(\zeta_6 z_1, \zeta_6^5z_2)$, (where
$\zeta_6$ is a 6-th primitive root of unity), and
$\beta:(z_1,z_2)\mapsto(-z_2,z_1)$. We observe that there are the
relations $\alpha_6^3=\beta^2$, $\alpha_6^6=\beta^4=1$,
$\alpha_6^{-1}\beta\alpha_6=\beta^{-1}$ so $\alpha_6$ and $\beta$ generate $\mathbb{D}_{12}\subset \Aut(A)$.
The points of $A$ with non trivial stabilizer for
$\mathbb{D}_{12}$ are the following:

\begin{eqnarray}\label{table: non trivial stab for D12}\end{eqnarray}
\vspace{-1cm}
$$
\begin{array}{c}
\mbox{ Points of $A$ with non trivial stabilizer for
$\mathbb{D}_{12}$}\\
\hline
\begin{array}{|c|rrl|}points\ in\ the\ same\ orbit & stabilizer&&\\
\hline
(0,0)&\mathbb{D}_{12}&=&\langle \alpha_6,\beta\rangle\\
\hline (0,\frac{1-\zeta_3}{3}), (\frac{-1+\zeta_3}{3},0),
(\frac{1-\zeta_3}{3},0), (0,\frac{-1+\zeta_3}{3})&\Z/3\Z&=&\langle
\alpha_6^2\rangle\\
\hline

(\frac{1-\zeta_3}{3},\frac{1-\zeta_3}{3}),
(\frac{-1+\zeta_3}{3},\frac{1-\zeta_3}{3}),(\frac{1-\zeta_3}{3},\frac{-1+\zeta_3}{3}),
(\frac{-1+\zeta_3}{3},\frac{-1+\zeta_3}{3})&\Z/3\Z&=&\langle
\alpha_6^2\rangle\\
\hline

(\frac{1}{2},\frac{1}{2}),
(\frac{\zeta_3}{2},\frac{1+\zeta_3}{2}),
(\frac{1+\zeta_3}{2},\frac{\zeta_3}{2})&\Z/4\Z&=&\langle\beta\rangle\\
\hline

(\frac{\zeta_3}{2},\frac{\zeta_3}{2}),
(\frac{1+\zeta_3}{2},\frac{1}{2}),
(\frac{1}{2},\frac{1+\zeta_3}{2})&\Z/4\Z&=&\langle\beta\rangle\\

\hline

(\frac{1+\zeta_3}{2},\frac{1+\zeta_3}{2}),
(\frac{1}{2},\frac{\zeta_3}{2}),
(\frac{\zeta_3}{2},\frac{1}{2})&\Z/4\Z&=&\langle\beta\rangle\\
\hline
(0,\frac{1}{2}), (0,\frac{\zeta_3}{2}),
(0,\frac{1+\zeta_3}{2}), (\frac{1}{2},0), (\frac{\zeta_3}{2},0),
(\frac{1+\zeta_3}{2},0)&\Z/2\Z&=&\langle\alpha_6^3\rangle\\

\hline

\end{array}\end{array}
$$

It follows that $F_{\mathbb{D}_{12}}$ is isometric to $D_5\oplus
A_2^2\oplus A_3^3\oplus A_1$.

First we consider the quotient by $\langle\alpha_6^2\rangle$,
which is a normal subgroup of $\mathbb{D}_{12}$. The quotient
$A/\langle \alpha_6^2\rangle$ is a surface with 9 singularities of type $A_2$, in the image of the points $p$ contained in the set
$\mathcal{P}:=\{(0,0)$, $(\frac{1-\zeta_3}{3},\frac{1-\zeta_3}{3})$, $
(\frac{-1+\zeta_3}{3},\frac{1-\zeta_3}{3})$, $(\frac{1-\zeta_3}{3},\frac{-1+\zeta_3}{3})$, $
(\frac{-1+\zeta_3}{3},\frac{-1+\zeta_3}{3})$, $
(\frac{1}{2},\frac{1}{2})$, $
(\frac{\zeta_3}{2},\frac{1+\zeta_3}{2})$, $
(\frac{1+\zeta_3}{2},\frac{\zeta_3}{2})\}$. This introduces 18
curves on the K3 surface $\widetilde{A/\langle\alpha_6^2\rangle}$,
desingularization of $A/\langle\alpha_6^2\rangle$: the curves
$K_{p}^j$, $p\in \mathcal{P}$, $j=1,2$ which desingularize the
point $p\in \mathcal{P}$. The automorphism $\beta\in \Aut(A)$
induces an automorphism $\beta'$ on
$\widetilde{A/\langle\alpha_6^2\rangle}$. Since $\beta$ fixes
$(0,0)$, $\beta'$ preserves the set of curves $\{K_{(0,0)}^{(j)}\}$,
$j=1,2$. The automorphism $\beta'$ fixes the intersection point
$K_{(0,0)}^{(1)}\cap K_{(0,0)}^{(2)}$ and switches the curves $K_{(0,0)}^{(1)}$
and $K_{(0,0)}^{(2)}$. The square $(\beta')^2$ preserves the curves
$K_{(0,0)}^{(1)}$ and $K_{(0,0)}^{(2)}$ and fixes their intersection point and another point on each curve. The points in
$\mathcal{P}-\{(0,0)\}$ have a trivial stabilizer with respect to
the action of $\langle\beta\rangle$ on $A$, so the 8 $A_2$-configurations generated by $K_{p}^j$, $p\in \mathcal{P}-\{(0,0)\}$, $j=1,2$, are
moved by $\beta'$. In particular neither $\beta'$ or $(\beta')^2$
have fixed points on these curves. The automorphism $\beta$ fixes
other 9 points of $A$ (see Table \eqref{table: non trivial stab for  D12}), which correspond to 3 points on
$\widetilde{A/\langle\alpha_6^2\rangle}-\{K_{p}^j\}$ (where $p\in
\mathcal{P}$, $j=1,2$) and thus to three singularities of type $A_3$ on $\widetilde{A/\langle\alpha_6^2\rangle}/\beta'$. The automorphism $\beta^2$ fixes other 6 points on $A$ (see Table \eqref{table: non trivial stab for  D12}), which correspond to 2 points of $\widetilde{A/\langle\alpha_6^2\rangle}-\{K_{p}^j\}$ 
(where $p\in\mathcal{P}$, $j=1,2$) and thus to one singular point of type $A_1$ on $\widetilde{A/\langle\alpha_6^2\rangle}/\beta'$.

Hence in the desingularization of $(\widetilde{A/\langle\alpha_6^2\rangle})/\beta'$, which is isomorphic to $X_{\mathbb{D}_{12}}$, there are the following
curves:\\ $K_{(0,0)}^{h}$, $h=0,\ldots,4$, which form a $D_5$; the
curves  $K_{(0,\frac{1-\zeta_3}{3})}^j$ and
$K_{(\frac{1-\zeta_3}{3},\frac{1-\zeta_3}{3})}^j$, $j=1,2$ which
form two disjoint copies of $A_2$ and which are image of the 8
copies of $A_2$ not preserved by $\beta'$; the curves $K_{(\frac{1}{2},\frac{1}{2})}^j$,
$K_{(\frac{\zeta_3}{2},\frac{\zeta_3}{2})}^j$, $K_{(\frac{1+\zeta_3}{2},\frac{1+\zeta_3}{2})}^j$, $j=1,2,3$,
which form 3 disjoint copies of $A_3$, and the curve $K_{(0,\frac{1}{2})}$ which is a copy of $A_1$.

\begin{figure}\scalebox{0.5}{
\begin{tikzpicture}
\node at (-1,1.5) {$K_{(0,0)}^{(4)}$};
\node at (0,1.5) {$K_{(0,0)}^{(0)}$};
\node at (1.5,1.5) {$K_{(0,0)}^{(1)}$};
\node at (2.5,1.5) {$K_{(0,0)}^{(3)}$};
\node at (3.5,0) {$K_{(0,0)}^{(2)}$};
    \draw (0,1) -- (0,-0.5);
    \draw (1.25,-0.25)--(-0.5,1);
    \draw (0.5,0)--(3,0);
    \draw (1.5,1)--(1.5,-0.5);
    \draw (2.5,1)--(2.5,-0.5);

\node at (5.5,-1) {$K_{(0,\frac{1-\zeta_3}{3})}^{(2)}$};
\node at (4.5,1.5) {$K_{(0,\frac{1-\zeta_3}{3})}^{(1)}$};
    \draw (4.5,-0.5)--(5.5,1);
    \draw (5.5,-0.5)--(4.5,1);

\node at (7.5,-1) {$K_{(\frac{1-\zeta_3}{3},\frac{1-\zeta_3}{3})}^{(2)}$};
\node at (6.5,1.5) {$K_{(\frac{1-\zeta_3}{3},\frac{1-\zeta_3}{3})}^{(1)}$};
    \draw (6.5,-0.5)--(7.5,1);
    \draw (7.5,-0.5)--(6.5,1);

\node at (10,1.5) {$K_{(\frac{1}{2},\frac{1}{2})}^{(1)}$};
\node at (11,-1) {$K_{(\frac{1}{2},\frac{1}{2})}^{(3)}$};
\node at (8.5,0) {$K_{(\frac{1}{2},\frac{1}{2})}^{(2)}$};
    \draw (9,0)--(12,0);
    \draw (10,-0.5)--(10,1);
    \draw (11,-0.5)--(11,1);

\node at (14.5,1.5) {$K_{(\frac{\zeta_3}{2},\frac{\zeta_3}{2})}^{(1)}$};
\node at (15.5,-1) {$K_{(\frac{\zeta_3}{2},\frac{\zeta_3}{2})}^{(3)}$};
\node at (13,0) {$K_{(\frac{\zeta_3}{2},\frac{\zeta_3}{2})}^{(2)}$};
    \draw (13.5,0)--(16.5,0);
    \draw (14.5,-0.5)--(14.5,1);
    \draw (15.5,-0.5)--(15.5,1);

\node at (19.5,1.5) {$K_{(\frac{1+\zeta_3}{2},\frac{1+\zeta_3}{2})}^{(1)}$};
\node at (20.5,-1) {$K_{(\frac{1+\zeta_3}{2},\frac{1+\zeta_3}{2})}^{(3)}$};
\node at (18,0) {$K_{(\frac{1+\zeta_3}{2},\frac{1+\zeta_3}{2})}^{(2)}$};
    \draw (18.5,0)--(21.5,0);
    \draw (19.5,-0.5)--(19.5,1);
    \draw (20.5,-0.5)--(20.5,1);

\node at (23,0) {$K_{(0,\frac{1}{2})}$};
    \draw (22.5,-0.5)--(22.5,1);

\end{tikzpicture}
}
\caption{Curves of $F_{\mathbb{D}_{12}}$ on
$X_{\mathbb{D}_{12}}$}\label{figure: A/D12 curves}
\end{figure}
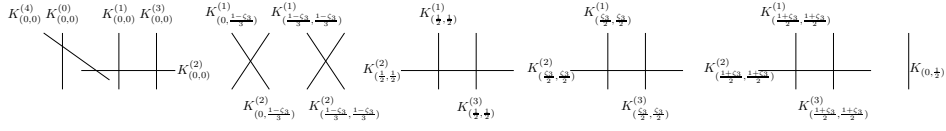

The intersection properties of these curves are presented in Figure \ref{figure: A/D12 curves}

The discriminant group of $\mathbb{D}_{12}$ is $(\Z/12\Z)^2\times (\Z/4\Z)^2\times \Z/2\Z$ generated by:
$d_1:=\frac{1}{4}\left(2K_{(0,0)}^{(4)}+K_{(0,0)}^{(1)}+2K_{(0,0)}^{(2)}+3K_{(0,0)}^{(3)}\right)$;\\
$d_2:=\frac{1}{4}\left(K_{(\frac{1}{2},\frac{1}{2})}^{(1)}+2K_{(\frac{1}{2},\frac{1}{2})}^{(2)}+3K_{(\frac{1}{2},\frac{1}{2})}^{(3)}\right)+\frac{1}{3}\left(K_{(0,\frac{1-\zeta_3}{3})}^{(1)}+2K_{(0,\frac{1-\zeta_3}{3})}^2\right)$;\\
$d_3:=\frac{1}{4}\left(K_{(\frac{\zeta_3}{2},\frac{\zeta_3}{2})}^{(1)}+2K_{(\frac{\zeta_3}{2},\frac{\zeta_3}{2})}^{(2)}+3K_{(\frac{\zeta_3}{2},\frac{\zeta_3}{2})}^{(3)}\right)+\frac{1}{3}\left(K_{(\frac{1-\zeta_3}{3},\frac{1-\zeta_3}{3})}^{(1)}+2K_{(\frac{1-\zeta_3}{3},\frac{1-\zeta_3}{3})}^2\right)$;\\
$d_4:=\frac{1}{4}\left(K_{(\frac{1+\zeta_3}{2},\frac{1+\zeta_3}{2})}^{(1)}+2K_{(\frac{1+\zeta_3}{2},\frac{1+\zeta_3}{2})}^{(2)}+3K_{(\frac{1+\zeta_3}{2},\frac{1+\zeta_3}{2})}^{(3)}\right)$;\\
$d_5:=\frac{1}{2}K_{(0,\frac{1}{2})}.$

The curves $K_{(0,0)}^{(j)}$, $j=1,2,3,4,5$, $K_{(\frac{1}{2},\frac{1}{2})}^{(h)}$, $K_{(\frac{\zeta_3}{2},\frac{\zeta_3}{2})}^{(h)}$, $K_{(\frac{1+\zeta_3}{2},\frac{1+\zeta_3}{2})}^{(h)}$, $h=1,2,3$ and $K_{(0,\frac{1}{2})}$
arise from the desingularization of the $\Z/4\Z$ quotient $\widetilde{A/\langle\alpha_6^2\rangle}\ra (\widetilde{A/\langle\alpha_6^2\rangle})/\beta'$
and so the class $v_{\mathbb{D}_{12}}:=d_1+9d_2+9d_3+d_4+d_5$ is contained in $NS(X_{\mathbb{D}_{12}})$, because $X_{\mathbb{D}_{12}}$ is the resolution of $\widetilde{A/\langle\alpha_6^2\rangle})/\beta'$.

Let us denote by $L_{\mathbb{D}_{12}}$ the lattice generated by the curves of $F_{\mathbb{D}_{12}}$ and by $v_{\mathbb{D}_{12}}$.
The discriminant group of $L_{\mathbb{D}_{12}}$ is $(\Z/12\Z)^2\times \Z/2\Z$ and it is generated by the vectors $\delta_1:=d_1+9d_2$,
$\delta_2:= d_1+9d_3$, $\delta_3:=2d_1+d_5$. Either $K_{\mathbb{D}_{12}}$ coincides with $L_{\mathbb{D}_{12}}$ or it is an overlattice of finite index
of $L_{\mathbb{D}_{12}}$. In the latter case there would be a non trivial vector $w$ in the discriminant group of $L_{\mathbb{D}_{12}}$, which is contained in $K_{\mathbb{D}_{12}}$.
Either $w$ or a multiple of $w$ generates either $\Z/2\Z$ or $\Z/3\Z$ in the discriminant group. It is easy to check that there is no a vector $w$ as required since it should corresponds either to the sum of $n$, $n\leq 7$, disjoint rational curves divided by $2$ or to the sum of $m$  disjoint $A_2$-configurations divided by 3, with $m\leq 2$. By Proposition \ref{prop: classes divisible by 2 and 3} these two possibilities are not acceptable, so $K_{\mathbb{D}_{12}}$ coincides with $L_{\mathbb{D}_{12}}$.

This result is different from the one given in \cite{W}: In our construction the lattice of Kummer type is generated by the classes of the curves arising from the desingularization of $A/G$ and by a class 4-divisible (i.e. the vector $v_{\mathbb{D}_{12}}$). In \cite[Proposition 2.1]{W} the lattice of Kummer type ($\Pi_{12}$, with the notation used in \cite{W}), is generated by the classes of the curves arising from the desingularization of $A/G$ and by a class 2-divisible (and not 4-divisible).  The discriminant group of the lattice $\Pi_{\mathbb{D}_{12}}$ described in \cite{W} is $(\Z/12\Z)^2\times (\Z/2\Z)^3$. This group has 5 generators. Since the rank of $\Pi_{\mathbb{D}_{12}}$ is 19 and the rank of $\Lambda_{K3}$ is 22, this is impossible because of Proposition \ref{prop: length of lattices}.

\subsection{The Kummer type lattices}\label{subsec: the Kummer type lattices}

Here we collect the results obtained above and the known ones in order to give a description of all the lattices of Kummer type. In particular we show that for all the lattices $K_G$ of Kummer type, the roots of $K_G$ coincide with the roots of $F_G$, which will be very useful in the following.

\begin{proposition}\label{prop: definition lattices FG and
KG}{\rm (See \cite{NikKummer} for $G=\Z/2\Z$; \cite{B} for
$G=\Z/n\Z$, $n=3,4,6$; Section \ref{subsection: non cyclic quotients of A} and \cite{W} for $G=\mathbb{D}_8$, $\mathbb{D}_8'$, $\mathbb{D}_{12}$, $\mathbb{T}$).} 

Let $A$ be an Abelian surface with an action of a finite group $G$ which does not contain translations. Let $X_G$ be the desingularization of $A/G$. If $X_G$ is a K3 surface, then $G$ is
one of the following 7 groups: $\Z/2\Z$, $\Z/3\Z$, $\Z/4\Z$,
$\Z/6\Z$,
$\mathbb{D}_8$, $\mathbb{D}_{12}$, $\mathbb{T}$. We recall that there are 2 different actions of the quaternion group denoted by $\mathbb{D}_8$ and $\mathbb{D}_8'$.

Let us assume that $X_G$ is a K3 surface (so $G$ is one of the 7
groups listed above). Let $K_i$ be the curves on $X_G$ arising by
the resolution of the singularities of $A/G$. Then the lattice
$F_G$ spanned by the curves $K_i$ is one of the following root
lattices:
\begin{eqnarray}
\begin{array}{l}
\begin{array}{|c||c|c|c|c|c|c|c|c|}
\hline
G&\Z/2\Z&/Z/3\Z&\Z/4\Z&\Z/6\Z\\
\hline F_G&A_1^{16}&A_2^9&A_3^4\oplus A_1^6&A_5\oplus A_2^4\oplus
A_1^5\\\hline\end{array}\\
$ $\\

\begin{array}{|c||c|c|c|c|c|}
\hline
G&\mathbb{D}_8&\mathbb{D}_8'&\mathbb{D}_{12}&\mathbb{T}\\
\hline F_G&D_4^2\oplus A_3^3\oplus A_1^2 & D_4^4\oplus
A_1^3&D_5\oplus A_3^3\oplus A_2^2\oplus A_1&E_6\oplus D_4\oplus
A_2^4\oplus A_1\\ \hline\end{array}
\end{array}
\end{eqnarray}

Let $K_G$ be the minimal primitive sublattice of $NS(X_G)$ which
contains the curves $K_i$, then $K_G$ is an overlattice of finite
index $r_G$ of $F_G$ with the following properties:
$$\begin{array}{|c||c|c|c|c|}
\hline
G&\Z/2\Z&\Z/3\Z&\Z/4\Z&\Z/6\Z\\
\hline
r_G&2^5&3^3&2^4&6\\
\hline
\rk(K_G)&16&18&18&18\\
\hline
K_G^{\vee}/K_G&(\Z/2\Z)^6&(\Z/3\Z)^3&(\Z/4\Z)^2\times(\Z/2\Z)^2&(\Z/6\Z)^4\\
\hline\end{array}$$ 
$$\begin{array}{|c||c|c|c|c|}
\hline
G&\mathbb{D}_8&\mathbb{D}_8'&\mathbb{D}_{12}&\mathbb{T}\\
\hline
r_G&2^3&2^4&4&3\\
\hline
\rk(K_G)&19&19&19&19\\
\hline
K_G^{\vee}/K_G&(\Z/4\Z)^3&(\Z/2\Z)^3&\Z/2\Z\times(\Z/12\Z)^2&(\Z/6\Z)^3\\
\hline
\end{array}$$ 
The roots of the lattice $K_G$
coincide with the roots of the lattice $F_G$ for all $G$.

By construction $K_G$ is a negative definite lattice primitively embedded in $NS(X_G)$ and
thus $\rho(X_G)\geq 1+\rk(K_G)$.
\end{proposition}
\proof The groups $G$ which act on $A$ in such a way that the reslution of $A/G$ is a K3 surface are classified by \cite[Lemma 3.3]{F}. The properties of $F_G$ and $K_G$ are proved in \cite[Section
1]{NikKummer} for $G=\Z/2\Z$; \cite[Section 1 and Theorem
2.5]{B} for $G=\Z/n\Z$, $n=3,4,6$; in Section \ref{subsection: non
cyclic quotients of A} and \cite[Proposition 2.1]{W} for $G=\mathbb{D}_8$, $\mathbb{D}_8'$, $\mathbb{D}_{12}$, $\mathbb{T}$. The unique
observation which has to be proved is that the root
system of $F_G$ coincides with the one of $K_G$. This was
explicitly proved in \cite[Proposition 1.3]{B} for
$G=\Z/n\Z$, $n=3,4,6$. In Section \ref{subsection: non cyclic quotients of A} we described a basis for $F_G$ and $K_G$ if $G$ is non cyclic and in \cite{NikKummer} a basis for $K_{\Z/2\Z}$ is given. One can explicitly write down a Gram matrix for lattice $K_G$. Since $K_G$ is a negative definite lattice, the number of vectors with a given self intersection is finite, and can be computed. In particular one computes the number of vectors of self-intersection $-2$ in $K_G$ (for example using the command \verb|ShortestVectors(-|$K_G$\verb|)| in \verb|Magma|) and one compares it with the number of vectors of self-intersection $-2$ in $F_G$. They coincide for every group $G$ in the list, and this concludes the proof.
\endproof

\subsection{The main results}\label{subsec: main on K3 covered by Abelian}

The aim of this section is to present and to prove our main result (Theorem \ref{theorem: X is A/G iff KG is in NS(X)}): one can deduce if a K3 surface is the quotient of an Abelian surface by checking if a certain lattice is primitively embedded in its N\'eron--Severi group. This essentially implies that one can construct the moduli space of the K3 surfaces which are desingularization of the quotient of an Abelian surface by a finite group as a moduli space of lattice polarized K3 surfaces. 

The other result of this section (Theorem \ref{theor: it suffices to have the curves on X}) is that one can deduce if a K3 surface is rationally $G$-covered by an Abelian surface by checking if a certain configuration of rational curves is present on the K3 surface.

We deduce by the combination of these two results a synthetic description of the lattices of Kummer type as overlattice with certain properties of the lattices $F_G$ (see Corollary \ref{cor: it suffices to have an overlattice of  F}).

\begin{theorem}\label{theorem: X is A/G iff KG is in NS(X)} Let $G$ be one of the following groups $\Z/n\Z$, $n=2,3,4,6$, $\mathbb{D}_{8}$, $\mathbb{D}_8'$, $\mathbb{D}_{12}$,
$\mathbb{T}$ and $K_G$ be the lattice of Kummer type defined above. A K3 surface is the minimal model of $A/G$ for a certain Abelian
surface $A$ if and only if $K_G$ is primitively embedded in
$NS(X_G)$.\end{theorem} \proof One of the implication is trivial:
if $X_G$ is the desingularization of $A/G$, then $NS(X_G)$
contains the classes of the curves arising from the
desingularization of $A/G$, so it contains the lattice $F_G$. By
definition $K_G$ is the minimal primitive sublattice of
$NS(X_G)$ which contains $F_G$ and so $K_G$ is primitively
embedded in $NS(X_G)$.

Let $X_G$ be a K3 surface such that $K_G$ is primitively embedded
in $NS(X_G)$. We first prove our result in case
$\rho(X_G)=1+\rk(K_G)$, i.e. it is the minimal possible. Let us
denote by $h$ the generator of the 1-dimensional subspace of
$NS(X_G)$ which is orthogonal to $K_G$, so $NS(X_G)$ is an
overlattice of finite index of $\Z h\oplus K_G$. Up to the action
of the Weyl group we can assume that $h$ is a pseudoample divisor
on $X_G$. Since $K_G$ is an overlattice of finite index of $F_G$,
$F_G$ is a root lattice and the roots of $F_G$ coincide with the
roots of $K_G$, the assumptions of Proposition \ref{prop:
the roots of L are smooth irreducible curves} (with $L:=K_G$ and
$R:=F_G$) are satisfied. Hence we can assume that the classes generating $F_G$
are supported on smooth irreducible rational curves. This fact
suffices to reconstruct the surface $A$ which is the minimal model
of the $G$-cover of $X_G$. This is well known in case $G=\Z/2\Z$,
see \cite{NikKummer}. The cases $G=\Z/3\Z$ and $G=\Z/4\Z$ are described in
\cite[Sections (4.1) and (4.2)]{B}. As example, we describe how
one can reconstruct $A$ in the cases $G=\Z/6\Z$,
$G=\mathbb{D}_8'$ and $G=\mathbb{T}$.

Let us assume $G=\Z/6\Z$. Then $F_G\simeq A_5\oplus A_2^4\oplus
A_1^5$ and $K_G$ is obtained by adding to $F_G$ the class
$$v:=\frac{1}{6}(\sum_{j=1}^5jK_1^{(j)})+\frac{1}{3}\sum_{i=2}^5(K_i^{(1)}+2K_i^{(2)})+\frac{1}{2}(\sum_{i=6}^{10} K_i^{(1)})$$
where $K_1^{(j)}$, $j=1,\ldots, 5$ is a basis of $A_5$,
$K_i^{(j)}$, $j=1,2$, $i=2,3,4,5$ is a basis of the $(i-1)$-th
copy of $A_2$, and $K_i^{(1)}$, $i=6,\ldots, 10$ is a generator of
the $(i-5)$-the copy of $A_1$. Let us now consider $3v$. It
exhibits the set of curves $\{K_1^{(1)}, K_1^{(3)}, K_1^{(5)},
K_6^{(1)}, K_7^{(1)},K_8^{(1)},K_9^{(1)},K_{10}^{(1)}\}$ as a set
of eight disjoint rational curves divisible by 2 on a K3 surface.
Then there exists a $2:1$ cover of $X_G$, $\widetilde{Y}\ra X_G$,
branched along these curves and such that the minimal model $Y$ of
$\widetilde{Y}$ is a K3 surface. The minimal model $Y$ is obtained
contracting the 8 $(-1)$-curves which are the $2:1$ cover of the
branch curves. Let us consider
the rational $2:1$ maps $\pi:Y\dashrightarrow X_G$. Then:
$\pi^{(-1)}(K_i^{(j)})$ splits in two rational curves for $j=1,2$, $i=2,3,4,5$, this gives 8 $A_2$-configurations on $Y$;
$\pi^{-1}(K_1^{(2)})$ is a rational curve which is a $2:1$ cover
of $K_1^{(2)}$ branched in two points; $\pi^{(-1)}(K_1^{(3)})$ is
a rational curve which is a $2:1$ cover of $K_1^{(3)}$ branched in
two points and we observe that after the contraction
$\widetilde{Y}\ra Y$, $\pi^{-1}(K_1^{(2)})$ and
$\pi^{-1}(K_1^{(3)})$ form a copy of $A_2$. So we have 9 copies of
$A_2$ on $Y$. By Proposition \ref{prop: known results on covers of K3 with  curves}, there exists an Abelian surface $A$ which is a 3:1 rational cover of $Y$. The minimal model of this
cover is an Abelian surface $A$, which is indeed a (rational) $G$-cover of $X_G$.

In case $\rho(X)>1+\rk(K_G)$, the proof follows by a standard deformation argument that we summarize here: If $\rho(X)>1+\rk(K_G)$, then $\rho(X)=20$. There exists a 1-dimensional family $\{X_t\}_{t\in \C}$, which deforms $X$, such that the generic member $X_t$ has Picard number 19 and $K_G$ is primitively embedded in $NS(X_t)$. Since generically $\rho(X_t)=1+\rk(K_G)$, for a generic $t$ there exists an Abelian surface $A_t$ which is a (rational) $G$-cover of $X_t$. In this way we produce a 1-dimensional family of Abelian surfaces $A_t$ admitting a $G$ action and such that the desingularization of $A_t/G$ is $X_t$. Generically $\rho(A_t)=3$ but there are special members, $A_{\overline{t}}$ in the family $\{A_t\}$ such that $\rho(A_{\overline{t}})=4$. These Abelian surfaces are (rational) $G$-covers of K3 surfaces $X_{\overline{t}}$, which has Picard number $20$. In particular there exist an Abelian surface $A$, special member of the family $\{A_t\}$, which is a (rational) $G$-cover of $X$.

Let us now consider the non-Abelian case. We remark that in this case $\rho(X)$ is necessarily equal to $1+\rho(K_G)$, because the latter is 20.
 
Let $G=\mathbb{D}_8'$ and $K_{\mathbb{D}_8}'$ as described in Section \ref{subsection: non cyclic quotients of A} and let $X$ be a K3 surface such that $K_{\mathbb{D}_8'}$ is primitively embedded in $NS(X)$. The classes $v_1'$ and $v_2'$ given in \eqref{eq: v1' and v2'} allows one to construct a $(\Z/2\Z)^2$-cover of $X$. Let us denote by $Y$ the minimal model of the $(\Z/2\Z)^2$-cover of $X$ branched with multiplicity 2 along each the 12 curves in the set $\mathcal{S}:=\{K_{(0,0)}^{(j)}$, $K_{(\frac{1+i}{2},0)}^{(j)}$, $K_{(\frac{i}{2},\frac{i}{2})}^{(j)}$, $K_{(\frac{1}{2},\frac{i}{2})}^{(j)}\}$, $j=1,2,3$. Let $\pi_Y:Y\dashrightarrow X$ be the rational map induced by the $(\Z/2\Z)^2$-cover and let $\langle \mu,\nu\rangle$ the cover group. We observe that $\pi_Y^{(-1)}(K_{(0,0)}^{(0)})$ consists of a unique irreducible rational curve and coincides with the inverse image of the $D_4$-configuration $K_{(0,0)}^{(j)}$, $j=0,1,2,3$. We denote this curve on $Y$ by $K_{(0,0)}$. Similarly the inverse images of the $D_4$-configuration $K_{(\frac{1+i}{2},0)}^{(j)}$, (resp. $K_{(\frac{i}{2},\frac{i}{2})}^{(j)}$,
$ K_{(\frac{1}{2},\frac{i}{2})}^{(j)}$) for  $j=0,1,2,3$, consists of a unique irreducible rational curve denoted by $K_{(\frac{1+i}{2},0)}$, (resp. $K_{(\frac{i}{2},\frac{i}{2})}$, $ K_{(\frac{1}{2},\frac{i}{2})}$). Since the curve $K_{(\frac{1}{2},0)}$ (resp. $K_{(\frac{1+i}{4},\frac{1+i}{4})}$, $K_{(\frac{1+i}{4},\frac{i-1}{4})}$) is not in the branch locus of the $(\Z/2\Z)^2$-cover of $X$, and does not meet the branch locus, its inverse image on $Y$ consists of 4 disjoint rational curves, denoted by $K_{(\frac{1}{2},0)}^{(j)}$ (resp. $K_{(\frac{1+i}{4},\frac{1+i}{4})}^{(j)}$, $K_{(\frac{1+i}{4},\frac{i-1}{4})}^{(j)}$ ) for $j=1,2,3,4$.
Thus on $Y$ there are 16 disjoint rational curves. Hence, by Proposition \ref{prop: known results on covers of K3 with  curves}, $Y$ is a Kummer surface of an Abelian surface $B$ and  there exists the following rational map
$\pi_B:B\dashrightarrow Y$, whose cover involution will be denoted by $\iota_B$. Hence there is a $8:1$ map, $\pi_Y\circ \pi_B:B\dashrightarrow X$. By construction the automorphisms $\mu$ and $\nu$ of $Y$ preserve the branch locus of the map $\pi_B:B\dashrightarrow Y$ and thus they induce two automorphisms $\mu_B$ and $\nu_B$ on $B$. Let us denote by $H_B$ the group generated by $\iota_B$, $\mu_B$ and $\nu_B$. By construction $B\dashrightarrow X$ is the map induced by the desingularization of the quotient $B/H_B$. In particular the group $H_B$ has order 8.

Let $\gamma:X\ra S$ be the contraction of all the curves in $F_{\mathbb{D}_8'}$. The singular surface $S$ has 4 singularities of type $D_4$ and 3 singularities of type $A_1$. It is immediate to check by our construction that $B\dashrightarrow X\stackrel{\gamma}{\ra} S$ coincides with the quotient $B\ra B/H_B$ and so $S=B/H_B$. The quotient singularities of type $D_4$ correspond to points whose stabilizer is the quaternion group, so the quaternion group $\mathbb{D}_8'$ has to be a subgroup of the group $H_B$, but the order of $H_B$ is 8, as the order of the quaternion group, so $H_B$ is the quaternion group. This implies that $X$ is the desingularization of the quotient $B/\mathbb{D}_8'$.

The case $G=\mathbb{D}_8$ is analogous. In case $G=\mathbb{D}_{12}$ one first considers a 4:1 cover of the K3 surface $X_G$. The minimal model of such a cover, say $Y$, contains 9 disjoint $A_2$-configurations, hence there exists an Abelian surface $B$ which is a 3:1 cover of $Y$. Then one proves that $X$ is the desingularization of $B/H_B$ where $H_B$ is a group generated by certain automorphism and, considering the singularities, one proves that $H_B$ must be $\mathbb{D}_{12}$ (since it has order 12 and has to contain $\mathbb{D}_{12}$).

Let us now consider the case $G=\mathbb{T}$. Let $X$ be a K3 surface such that $K_\mathbb{T}$ is primitively embedded in $NS(X)$. So there are 19 curves which span the lattice $E_6\oplus D_4\oplus A_1\oplus A_2^4$ and there is a 3 divisible class which involves 6 disjoint $A_2$-configurations. So there exists a $3:1$ cover of $X$ whose minimal model is a K3 surface $Y$. We denote by $\pi:Y\dashrightarrow X$ the $3:1$ rational map. The inverse image on $Y$ of the curves in $F_{\mathbb{D}_{12}}$ consists of 19 rational curves which span the lattice $D_4^4\oplus A_1^3$ (where a copy of $D_4$ is mapped by $\pi$ to the $E_6$ contained in $F_{\mathbb{T}}$, three other copies of $D_4$ are mapped to the unique copy of $D_4$ in $F_{\mathbb{T}}$, the three copies of $A_1$ are mapped by $\pi$ to the unique copy of $A_1$ in $F_{\mathbb{T}}$). We observe that $D_4^4\oplus A_1^3\simeq F_{\mathbb{D}_8'}$. In order to reconstruct the Abelian surface which is the cover of $X$ it suffices to prove that not only $F_{\mathbb{D}_8'}$, but exactly $K_{\mathbb{D}_8'}$ is primitively embedded in $NS(Y)$. Once one proves this, one finds an Abelian surface $B$ such that $Y$ is the minimal resolution of $B/\mathbb{D}_8'$ (we already proved this result) and one deduces that $X$ is the minimal resolution of $B/\mathbb{T}$ as in the previous cases. In Section \ref{subsection: non cyclic quotients of A} we constructed the lattice $K_{\mathbb{D}_8'}$ introducing four divisible vectors. Two of them ($v_1'$ and $v_2'$) are strictly related with the geometry of the quotient that we are considering. The property of these two vectors, which is essential in order to reconstruct the Abelian surface $B$ with a $\mathbb{D}_8'$-action, is that the curves appearing in these two divisible classes are all contained in the $D_4$-configurations, i.e. the curves which generates the three copies of $A_1$ in $F_{\mathbb{D}_8'}$ do not appear in these divisible classes. Since there are 15 disjoint rational curves contained in the set of the 19 curves which span $D_4^4\oplus A_1^3$, we know that there are also 4 independent divisible classes in $NS(Y)$, by Proposition \ref{prop: known results on covers of K3 with  curves}. Now we have to show that at least 2 of them can be chosen to have no components in the direct summands $A_1^3$ of $D_4^4\oplus A_1^3$. Suppose the opposite, this means that there is a choice of three  divisible vectors $n_1$, $n_2$ and $n_3$ such that all the elements in $\langle n_1,n_2,n_3\rangle$ have components among the generators of $A_1^3$. Just to fix the notation we gives to the curves in $F_{\mathbb{D}_8'}$ the same name as in Section \ref{subsection: non cyclic quotients of A}. We chose the first class $n_1$ in such a way that it has some components in $A_1^3$. We recall that the divisible classes are the sum of eight disjoint rational curves divided by 2 and that they are linear combinations of the elements of the discriminant group. We observe that a divisible class has components among the generators of $A_1^3$ if and only if at least one of the vectors (of the discriminant group) $d_9$, $d_{10}$ and $d_{11}$ appears with a non trivial coefficient in its expression. Since the generators of the discriminant group $d_i$  with $i\neq 9,10,11$ are the sum of two rational curves divided by 2, in the expression of $n_1$ an even number of vectors $d_9$, $d_{10}$ and $d_{11}$ with a non trivial coefficient appears. So we can assume that $n_1:=d_9+d_{10}+m_1$, where $m_1\in \langle d_j\rangle$, $j=1,\ldots 8$. Now we construct a second divisible class $n_2$ assuming that it has some components among the generators of  $A_1^3$. If $n_2:=d_9+d_{10}+m_2$, where $m_2\in \langle d_j\rangle$, $j=1,\ldots 8$, then $n_1+n_2\in \langle d_j\rangle$, $j=1,\ldots 8$, i.e. it has no components among the generators of $A_1^3$. So we can assume that $n_2:=d_9+d_{11}+m_2$, $m_2 \in \langle d_j\rangle$, $j=1,\ldots 8$. We observe that $n_1+n_2=d_{10}+d_{11}+m_3$, $m_3 \in \langle d_j\rangle$, $j=1,\ldots 8$.
But now there is no way to choose $n_3$ is such a way that all the elements in $\langle n_1,n_2,n_3\rangle$ have components among the curves generating $A_1^3$. Indeed every pair of elements in $\{d_9,d_{10},d_{11}\}$ appears with a non trivial coefficient in $n_1$ or in $n_2$ or in $n_1+n_2$. 
This proves that if on a K3 surface $Y$ there is configuration of 19 rational curves which span the lattice $F_{\mathbb{D}_8'}$, then the lattice $K_{\mathbb{D}_8'}$ is primitively embedded in $NS(Y)$ and so $Y$ is the minimal resolution of the quotient of an Abelian surface $B$ by the group $\mathbb{D}_8'$. Moreover, this concludes the proof in the unique remaining case $G=\mathbb{T}$.

\endproof

\begin{rem} In \cite{B} the proof of the previous result is given in case $G$ is a cyclic group of order greater then
2. The proof given in case $\rho(X_G)$ is the minimal possible
coincides with our proof. In case $\rho(X_G)$ is greater (and
indeed 20), in \cite{B} it is observed that one can use a deformation argument as we did, but an alternative proof is given.
Unfortunately, it is based on \cite[Lemma 3.2]{B}, which contains
a mistake. Indeed, using the notation of \cite[Lemma 3.2]{B}, it
is true that there exists an orthogonal embedding $\eta$ of
$\{A_{k_1},\ldots A_{k_n}\}$ in a system of roots, $Q$, of type
$\mathbb{A}$ such that (up to the action of the Weyl group),
$\eta(A_{k_i})$ is contained in a chosen basis of $Q$ for
every $i=1,\ldots n$, but the same result is not necessarily true
if the system of roots $Q$ is of type $\mathbb{D}$. A simple
counterexample is given by the orthogonal embedding of
$\{A_1,A_1,A_1,A_1\}$ in $D_4$ given by $\{\epsilon_1+\epsilon_2,
\epsilon_1-\epsilon_2, \epsilon_3+\epsilon_4,
\epsilon_3-\epsilon_4\}$ which can not be contained in a basis of
$D_4$ (up to the action of the Weyl group of $D_4$).
\end{rem}

The advantage of the result in Theorem \ref{theorem: X is A/G iff KG is in NS(X)} is that one relates a purely geometric property with a purely lattice theoretic property. This is what is needed in order to describe the lattice polarized moduli space of the K3 surfaces with a certain geometric property, so we immediately obtain the following corollary:
\begin{corollary}\label{corollary: L polarized A/G}
Let $\mathcal{L}_G$ be the set of lattices $L_G$ satisfying:\begin{itemize}\item[(a)] $L_G$ has rank $1+\rk(K_{G})$,\item[(b)] $L_G$ is hyperbolic,\item[(c)] $L_G$ admits a primitive embedding in $\Lambda_{K3}$, \item[(d)] there exists a primitive embedding of $K_G$ in $L_G$.\end{itemize} A K3 surface is the desingularization of the quotient of an Abelian surface by $G$ if and only if it is an $L_G$-polarized K3 surface for an $L_G\in\mathcal{L}_G$. 

In particular the coarse moduli space of the K3 surfaces which are desingularization of the quotient $A/G$ for an Abelian surface $A$ has infinitely many components of dimension $19-\rk(K_G)$.\end{corollary}

We observe that the conditions $(a)$, $(b)$ and $(d)$ in Corollary \ref{corollary: L polarized A/G} imply that $L_G$ is an overlattices of finite index $l_G$ of $\Z h\oplus K_{G}$, where $h$ is a vector with a positive self intersection $h^2$. The condition $(c)$ implies that $h^2$ is even and imposes several restriction to $l_G$. The concrete possibilities for the lattices in $\mathcal{L}_G$ are classically known for $G=\Z/2\Z$ (see for example \cite[Theorem 2.7]{GS} for a recent reference) and for $G=\Z/3\Z$, see \cite{Barthclass}.\\

In \cite{NikKummer}, it is proved that it is not necessary to check the existence of a primitive embedding of $K_{\Z/2\Z}$ in the N\'eron--Severi group of a K3 surface to conclude that it is an Abelian surface: it suffices to know that it contains 16 disjoint smooth irreducible rational curves. We underline that from the point of view of the description of the moduli space this result is not very useful, because we have no a way to translate the condition "certain $-2$ classes correspond to irreducible curves" in the context of the lattice polarized K3 surfaces. On the other hand this results is very nice from a geometric point of view, since it can be stated also in the following way: if a K3 surface admits a model with 16 nodes, then it is a Kummer surface (for example this is can be used to conclude that a quartic with 16 nodes is a Kummer surface). The similar result was generalized to the group $G=\Z/3\Z$ by Barth in \cite{Barthnine}. Here we generalize this result to all the other admissible groups.

\begin{theorem}\label{theor: it suffices to have the curves on X}  Let $G$ be one of the following groups $\Z/n\Z$, $n=2,3,4,6$, $\mathbb{D}_8$, $\mathbb{D}_8'$, $\mathbb{D}_{12}$, $\mathbb{T}$
and $F_G$ be the lattice defined above. Then a K3 surface is the
minimal model of $A/G$ for a certain Abelian surface $A$ if and
only if $F_G$ is embedded in $NS(X_G)$ and there
exists a basis of $F_G$ which represents irreducible smooth curves
on $X_G$.\end{theorem}
\proof This result is known if $G=\Z/2\Z$, see \cite{NikKummer} and if $G=\Z/3\Z$, see \cite{Barthnine}. In the proof of the Theorem 
\ref{theorem: X is A/G iff KG is in NS(X)} we proved the statement in case $G=\mathbb{D}_8'$. Here we give a complete proof in the case $G=\Z/4\Z$. The other cases are very similar. The lattice $F_{\Z/4\Z}$ has rank 18 and length 10. Since the length of a lattice of rank 18 primitively embedded in $\Lambda_{K3}$ is at most 4 (=22-18) we know that $F_{\Z/4\Z}$ is not primitively embedded in $\Lambda_{K3}$ and so there is an overlattice of finite index of $F_{\Z/4\Z}$, called $R_{\Z/4\Z}$, which is primitively embedded in $\Lambda_{K3}$. In order to construct an overlattice $R_{\Z/4\Z}$ of $F_{\Z/4\Z}$ we have to add to $F_{\Z/4\Z}$ certain elements which are non trivial in the discriminant group of $F_{\Z/4\Z}$ and which have an even self intersection. Moreover we have to recall that if the sum of $m$ disjoint rational curves is divided by 2, then $m$ is either 16 or 8. 

Let us consider the lattice $F_{\Z/4\Z}=A_3^4\oplus A_1^6$. We denote by $a_i^{(j)}$, $i=1,2,3$, $j=1,2,3,4$ the basis of the 
$j$-th copy of $A_3$ and by $a^{(j)}$, $j=5,6,7,8,9,10$ the generator of the $(j-4)$-th copy of $A_1$. The discriminant of $F_{\Z/4\Z}$ is generated by $d_j:=\frac{1}{4}\left( a_1^{(j)}+2a_2^{(j)}+3a_3^{(j)}\right)$, $j=1,2,3,4$, $d_j:=\frac{a^{(j)}}{2}$, $j=5,6,7,8,9,10$. 

Since $l(F_{\Z/4\Z})-l(R_{\Z/4\Z})$ has to be at least 6, we have to add at least 3 divisible vectors to $F_{\Z/4\Z}$ in order to obtain $R_{\Z/4\Z}$. First we suppose to add three vectors, $v_1,v_2,v_3$ such that $\langle v_1,v_2,v_3\rangle=(\Z/2\Z)^3$ in the discriminant group (i.e. no vectors among $v_1$, $v_2$ ,$v_3$ has order 4 in the discriminant group of $F_{\Z/4\Z}$). 
Every vector which generates $\Z/2\Z$ in the discriminant group of $F_{\Z/4\Z}$ is a linear combination of $2d_j$ for $j=1,2,3,4$ and $d_k$ for $k=5,\ldots,10$. The curves which appear with a non trivial coefficient in each of these linear combinations are among the 14 disjoint rational curves $\{a_1^{(j)}, a_3^{(j)}, a^{(k)}\}$ for $j=1,2,3,4$ and $k=5,\ldots,10$. We recall that it is possible to add three independent divisible 2-classes starting from 14 disjoint rational curves, but it is not possible to add 4 independent divisible classes using only 14 rational curves. So we can add exactly the 3 vectors $v_1$, $v_2$ and $v_3$. Up to permutations of the indices the unique possibility for the 3 vectors $v_1$, $v_2$ and $v_3$ is $v_1:=2(d_1+d_2+d_3+d_4)$, $v_2:=2d_1+2d_2+d_5+d_6+d_7+d_8$, $v_3:=2d_1+2d_3+d_7+d_8+d_9+d_{10}$. The lattice $R_{\Z/4\Z}$ obtained adding to $F_{\Z/4\Z}$ the vectors $v_1$, $v_2$ and $v_3$ is an overlattice of index $2^3$. One can directly compute its discriminant group and one finds that the discriminant group of this lattice is $(\Z/4\Z)^2\times (\Z/2\Z)^4$. But the length of this lattice is 6, which is not admissible. 

We conclude that there is at least one vector, say $v_1$ in $F_{\Z/4\Z}/R_{\Z/4\Z}$ which generates a copy of $\Z/4\Z$ in the discriminant group of $F_{\Z/4\Z}$. We recall that $(v_1)^2$ has to be an even number, that $(d_j)^2=-3/4$ if $j=1,2,3,4$ and that $(d_k)^2=-1/2$ if $k=5,\ldots, 10$. Moreover, $2v_1\mod F_{\Z/4\Z}$ has to be the sum of 8 disjoint rational curves divided by 2 (since the sum of $n$ rational curves can not divided by 2 if $n\leq 14$ and $n\neq 8$). So there are only the following 2 possibilities modulo $F_{\Z/4\Z}$ (up to a permutation of the indices): either\\
$(i)$ $v_1:=d_1+d_2+d_3+d_4+d_5+d_6+d_7+d_8+d_9+d_{10}$ or\\
$(ii)$ $v_1:=d_1+d_2+d_3+d_4+d_5+d_6$.\\
In case $(i)$ one can construct a $4:1$ cover of $X$ whose branch divisor is $v_1$. So we have a map $Y\ra X$ which is 4:1. By construction the minimal model of $Y$ has a trivial canonical bundle and its Euler characteristic is 0, so this surface is an Abelian surface and we conclude the proof. We remark that it suffices to observe that the divisor $v_1$ in case $(i)$ is the one described by Bertin in \cite[Page 270]{B} where it is proved that the minimal model of a $4:1$ cover of a K3 surface whose branch locus has a certain property has to be an Abelian surface. 

In case $(ii)$ the $4:1$ cover associated to the vector $v_1$ produces a K3 surface, and not an Abelian surface. Thus we have to analyze not only the vector $v_1$, but also the vectors $v_2$ and $v_3$ in order to show that $R_{\Z/4\Z}$ coincides with $K_{\Z/4\Z}$. We now consider the vectors $v_2$ and $v_3$. Up to replace, possibly, $v_2$ (resp. $v_3$) with $2v_2$ (resp. $2v_3$), we have that $v_2$ (resp. $v_3$) generates a copy of $\Z/2\Z$ and consists of the sum of 8 disjoint rational curves divided by 2; 4 of these curves have to be chosen among $\{a_1^{(j)},a_3^{(j)}\}$, $j=1,2,3,4$ since these are the eight disjoint rational curves of the divisible vector $2v_1$. Up to a permutation of the indices we can assume that $v_2:=2d_1+2d_2+d_5+d_6+d_7+d_8$ and $v_3:=2d_1+2d_3+d_7+d_8+d_9+d_{10}$. Now we consider the vector $v_1+v_3$ (which is surely contained in $R_{\Z/4\Z}$). It is $3d_1+d_2+3d_3+d_4+d_5+d_6+d_7+d_8+d_9+d_{10}$. Modulo $F_{\Z/4\Z}$ and a change of the index of the generators of $A_3$, this coincides with the vector $v_1$ in case $(i)$. So the minimal model of $4:1$ cover of $X$ whose branch divisor is $v_1+v_3$ is an Abelian surface and we conclude the proof as before. 

The other cases are similar (but easier): one checks that the length of $F_{G}$ is greater than $22-\rk(F_G)$, one deduces that one has to add some divisible classes in order to construct the lattice $R_{G}$ which is the minimal primitive sublattice of $\Lambda_{K3}$ containing $F_{G}$. One identifies these classes (recalling the condition that they are linear combinations of elements of the discriminant group of $F_{G}$ and the conditions imposed by Proposition \ref{prop: known results on covers of K3 with  curves}). Then one compares the lattice $R_G$ with $K_{G}$ or one explicitly construct a certain cover of $X$ in order to show either that $R_G=K_G$ (which implies that $X$ is the desingularization of $A/G$ by Theorem \ref{theorem: X is A/G iff KG is in NS(X)}) or directly that there exists an Abelian surface $A$, such that $X$ is the resolution of $A/G$.\endproof

\begin{corollary}\label{cor: it suffices to have an overlattice of F} Let $G$ be one of the groups $\Z/n\Z$, $n=2,3,4,6$, $\mathbb{D}_8$, $\mathbb{D}_8'$, $\mathbb{D}_{12}$, $\mathbb{T}$
and $F_G$ be the lattice defined above. Let $H_G$ be the minimal primitive sublattice of $\Lambda_{K3}$ which contains $F_G$ and such that the root lattice of $F_G$ coincides with the one of $H_G$. Then $H_G\simeq K_G$.\end{corollary}
\proof By hypothesis $H_G$ is a negative definite lattice primitively embedded in $\Lambda_{K3}$ and $\rk(H_G)=\rk(F_G)$. Let $D$ be a vector in $\Lambda_{K3}$ which is orthogonal to $H_G$ and has a positive self intersection. By the Torelli theorem there exists a K3 surface, $X$, whose transcendental lattice is the orthogonal to $\Z D\oplus H_G$ in $\Lambda_{K3}$. The N\'eron--Severi group of $X$ is an overlattice of finite index of $\Z D\oplus H_G$ such that $H_G$ is primitively embedded in it. Under our assumptions of $H_G$ we can apply Proposition \ref{prop: the roots of L are smooth  irreducible curves} to $L=H_G$ and $R=F_G$. So the lattice $F_G$ is spanned by irreducible rational curves on $X$. By Theorem \ref{theor: it suffices to have the curves on X}, it follows that $X$ is the desingularization of the quotient $A/G$ for a certain Abelian surface $A$. In this case the minimal primitive sublattice which contains the curves of the lattice $F_G$ is $K_G$, but by the hypothesis the minimal primitive sublattice of $NS(X)\subset \Lambda_{K3}$ which contains $F_G$ is $H_G$, so $K_G$ coincides with $H_G$.\endproof

\begin{rem}
The hypothesis that the roots of $H_G$ coincide with the ones of
$F_G$ in Corollary \ref{cor: it suffices to have an overlattice of  F} is essential. Indeed let us
consider the case $G=\Z/2\Z$. The lattice $F_G$ is $A_1^{16}$ and
let us denote by $K_i$, $i=1,\ldots 16$ the generators of this
lattice. Let us consider the vectors
$v_j:=(\sum_{i=1}^4K_{4j+i})/2$, $j=0,1,2,3$,
$w_1:=(K_1+K_2+K_5+K_6+K_9+K_{10}+K_{13}+K_{14})/2$ and
$w_2:=(K_1+K_3+K_5+K_7+K_9+K_{11}+K_{13}+K_{15})$. Let us denote by $H_{\Z/2\Z}$ the lattice obtained adding to $F_G$ the vectors $v_i$, $i=0,1,2,3$ and $w_h$, $h=1,2$. It is an overlattice (of index $2^6$) of $F_G$ which admits a primitive embedding in $\Lambda_{K3}$, but it is not isometric to the Kummer lattice (which in fact is an overlattice of index $2^5$ of $F_G$). In this case $v_1$ is a root of $H_G$ which is not contained in $F_G$.
\end{rem}

\subsection{K3 surfaces (rationally) $\Z/3\Z$-covered by Abelians surfaces}\label{K3 rationally covered by Z/3Z} $ $

In \cite{GS} it is observed that every Kummer surface $Km(A)$ (i.e. every K3 surface which is the desingularization of $A/\left(\Z/2\Z\right)$) admits the group $(\Z/2\Z)^4$ as group of symplectic automorphisms. Moreover, $Km(A)$ is also the quotient of a K3 surface by the symplectic action of $(\Z/2\Z)^4$. This result is based on the observation that if a K3 surface is a Kummer surface $Km(A)$, then the translations by the 2 torsion points of $A$ induce symplectic automorphisms on $Km(A)$.

A similar result can be obtained if the K3 surface $X_G$ is the (desingularization of the) quotient of an Abelian surface by an action of the group $\Z/3\Z$. 
\begin{proposition}\label{prop: A/Z3 admits Z/3^2} Let $X$ be the desingularization of the quotient of an Abelian surface $A$ by the group $\Z/3\Z$. Then $X$ admits a symplectic action of the group $(\Z/3\Z)^2$. Moreover, there exists a K3 surface $Y$ which admits a symplectic action of $(\Z/3\Z)^2$ such that $X$ is the desingularization of $Y/(\Z/3\Z)^2$.\end{proposition}
\proof Let $A$ be an Abelian surface admitting an automorphism $\alpha_A$ of order 3 such that $X$ is the desingularization of $A/\alpha_A$. Let $A[3]$ be the group of 3 torsion points of $A$ and let $\langle P, Q\rangle\subset A[3]$ be the set of points fixed by $\alpha_A$. Let us denote by $t_{P}$ and $t_{Q}$ the translation on $A$ by the points $P$ and $Q$ respectively. Then $(\Z/3\Z)^2\simeq \langle t_{P}, t_{Q}\rangle\subset\Aut(A)$ and the automorphisms $t_{P}$ and $t_{Q}$ commute with $\alpha_A$. So $t_{P}$ and $t_{Q}$ induce two automorphisms of order 3 on $A/\alpha_A$ which lifts to two automorphisms, $\tau_{P}$ and $\tau_{Q}$, on $X$. The period of $X$ (i.e. the generator of $H^{2,0}(X)$) is induced by the period of $A$, which is preserved by the translations. So $\tau_{P}$ and $\tau_{Q}$ are symplectic automorphisms of $X$. This gives a symplectic action of $(\Z/3\Z)^2$ on $X$.

On the other hand, $X$ contains 9 disjoint $A_2$-configurations of rational curves (which generates the lattice $F_{\Z/3\Z}$) and the minimal primitive sublattice $K_{\Z/3\Z}$ which contains all these curves contains also several divisible classes. In particular, let us denote by $a_i^{(j)}$, $i=1,2$, $j=1,\ldots,9$ the basis of the $j$-th copy of $A_2$. Up to a choice of the indices, $K_{\Z/3\Z}$ contains also the classes (mod $F_{\Z/3\Z}$) $$v_1:=\frac{1}{3}\left(\sum_{i=1}^6 a_1^{(j)}-a_2^{(j)}\right);\  v_2=\frac{1}{3}\left(\sum_{j=1}^2 (a_1^{(j)}-a_2^{(j)})-\sum_{h=3}^4(a_1^{(h)}-a_2^{(h)})+\sum_{k=7}^8 (a_1^{(k)}-a_2^{(k)})\right),$$
as shown in \cite[Page 269]{B} with a slightly different notation.
But the presence of these divisible classes allows one to reconstruct a $(\Z/3\Z)^2$ cover of $X$ (one first constructs the $3:1$ cover associated to the class $v_1$ as in Section \ref{subsec: lattice} and then one considers the pull back of the class $v_2$, which allows one to construct another $3:1$ cover). With this process one obtains a non minimal surface, whose minimal model $Y$ is a K3 surface which is a (rational) $(\Z/3\Z)^2$-cover of $X$, hence $X$ is the desingularization of the quotient of the K3 surface $Y$ by the group $(\Z/3\Z)^2$. \endproof

\begin{corollary}
The 1-dimensional families of K3 surfaces which are desingularizations of the quotients $A/\Z/3\Z$ for certain Abelian surfaces $A$ are contained in the intersection between the 3-dimensional families of the K3 surfaces which are (desingularization of) quotients of K3 surfaces by a symplectic action of $(\Z/3\Z)^2$ and the 3-dimensional families of K3 surfaces which admit a symplectic action of $(\Z/3\Z)^2$.
\end{corollary}

\begin{rem} The existence of the surface $Y$ in the Proposition \ref{prop: A/Z3 admits Z/3^2} directly follows by the primitive embedding of lattice $M_{(\Z/3\Z)^2}$ in the lattice $K_{\Z/3\Z}$, after proving Theorem \ref{theorem: S is A/G iff MG is in NS(S)}. Similarly one obtains that if $X$ is the minimal model of the quotient $A/(\Z/4\Z)$ for a certain Abelian surface then it is also the minimal model of the quotient $Y/(\Z/4\Z)$ for a certain K3 surface $Y$, since $M_{\Z/4\Z}\subset K_{\Z/4\Z}$. \end{rem}

In Proposition \ref{prop: A/Z3 admits Z/3^2} we proved that a K3 surface $X$ which is (rationally) $\Z/3\Z$-covered by an Abelian surface, necessarily admits certain symplectic automorphisms induced by translation on the Abelian surface. Here we observe that there exists another automorphism on $A$ which induces symplectic automorphism on $X$.

\begin{proposition}\label{prop: A/Z3 admits involution}
Let $X$ be a K3 surface such that $K_{\Z/3\Z}$ is primitively
embedded in $NS(X)$, then $X$ admits a symplectic involution
$\iota_X$ such that $K_{\Z/6\Z}$ is primitively embedded in
$NS(W)$ where $W$ is the K3 surface 
minimal model of $X/\iota_X$. \end{proposition}
\proof Every Abelian surface admits an involution $\iota_A:A\ra A$ which sends every point to its inverse with respect to the group law of $A$. Under the hypothesis on $X$ there exists an Abelian surface $A$ with an automorphism $\alpha_A\in\Aut(A)$ of order 3 such that $X$ is the desingularization of $A/\alpha_A$. The automorphism  $\iota_A$ and $\alpha_A$ commutes and generate an automorphism $\alpha_A\circ\iota_A$ of order 6 which preserves the period of $A$. The involution $\iota_A$ induces an involution $\iota_X$ on $X$. The singular surface  $A/(\alpha_A\circ \iota)$ is birational to $X/\iota_X$. Since the minimal model of $A/(\alpha_A\circ\iota)$ is a K3 surface, also the minimal model of $X/\iota_X$ is a K3 surface and these surfaces are isomorphic. We call this surface $W$ and we observe that it is constructed as minimal model of the quotient of an Abelian surface by the action of $\Z/6\Z=\langle \alpha_A\circ\iota\rangle$, so $K_{\Z/6\Z}$ is primitively embedded in $NS(W)$.\endproof

A generalization of the previous result can be done substituting $(\Z/3\Z, \Z/6\Z)$ with $(\mathbb{D}_8', \mathbb{T})$:
\begin{corollary}
Let $S$ be a K3 surface such that $K_{\mathbb{D}_8'}$ is primitively
embedded in $NS(S)$, then $S$ admits an automorphism of order 3, $\gamma_S$, such that $K_{\mathbb{T}}$ is primitively embedded in
$NS(\widetilde{S/\gamma_S})$ where $\widetilde{S}/\gamma_S$ is the
minimal resolution of $S/\gamma_S$. \end{corollary}

Putting together the Propositions \ref{prop: A/Z3 admits Z/3^2} and \ref{prop: A/Z3 admits involution}, one obtains the following corollary
\begin{corollary} Let $X$ be a K3 surface which is (rationally) $(\Z/3\Z)$-covered by an Abelian surface. The group $\mathfrak{A}_{3,3}$ acts symplectically on $X$.\end{corollary}
\proof It suffices to prove that the involution $\iota_A$ and the translations $t_{P}$ and $t_Q$ introduced in proofs of Propositions \ref{prop: A/Z3 admits involution} and \ref{prop: A/Z3 admits Z/3^2} generate $\mathfrak{A}_{3,3}$. This can be easily checked, for example one can specialize the Abelian surface $A$ to the product of two elliptic curves with $j$-invariant equal to 0. The order 3 automorphism $\alpha_A$ (defined in proof of Proposition \ref{prop: A/Z3 admits Z/3^2}) fixes the points $(0,0)$, $P:=(\frac{1}{3}(1-\zeta_3),0)$ and $Q:=(0,\frac{1}{3}(1-\zeta_3))$. This identifies the translation $t_P$ and $t_Q$ and it is immediate to verify that $\langle t_P,\iota\rangle\simeq \langle t_Q,\iota\rangle$ is the dihedral group of order 6 and then $\langle t_P,t_Q,\iota\rangle$ is $\mathfrak{A}_{(3,3)}$.\endproof

\section{K3 surface quotients of K3 surfaces}\label{sc: K3 covered by K3}
The aim of this section is to extend some of the results proved for the K3 surfaces which are (rationally) covered by an Abelian surface, to the K3 surfaces which are (rationally) covered by a K3 surface. We will denote by $Y_G$ a K3 surface which admits a symplectic action of the group $G$ and by $S_G$ the minimal resolution of the quotient $Y_G/G$. It is well known that $S_G$ is a K3 surface (see \cite{Niksympl}).

\begin{proposition}\label{prop: EG and MG}

Let $Y_G$ be a K3 surface and $G\in\Aut(Y_G)$ be a finite group. Let $S_G$ be the minimal model of $Y_G/G$. Then $S_G$ is a K3 surface if and only if $G$ acts symplectically on $Y_G$. If $G$ is Abelian, then it is one of the following 14 groups $\Z/n\Z$, $n=2,\ldots,8$, $(\Z/m\Z)^2$, $m=2,3,4$, $\Z/2\Z\times \Z/t\Z$, $t=4,6$, $(\Z/2\Z)^j$, $j=3,4$. 

Let $M_i$ be the curves on $S_G$ arising from
the resolution of the singularities of $Y_G/G$. Then the lattices
$E_G$ spanned by the curves $M_i$ is one of the following root
lattices:
\begin{eqnarray}
\begin{array}{l}
\begin{array}{|c||c|c|c|c|c|c|c|c|c|c|c|c|c|c|}
\hline
G&\Z/2\Z&/Z/3\Z&\Z/4\Z&\Z/5\Z&\Z/6\Z&\Z/7\Z&\Z/8\Z\\
\hline E_G&A_1^8&A_2^6&A_3^4\oplus A_1^2&A_4^4&A_5^2\oplus A_2^2\oplus A_1^2& A_6^3&
A_7^2\oplus A_3\oplus A_1\\
\hline\end{array}\\
$ $\\
\begin{array}{|c||c|c|c|c|c|c|c|c|c|c|c|c|c|c|c|c|c|}
\hline
G&(\Z/2\Z)^2&(\Z/2\Z)^3&(\Z/2\Z)^4&\Z/2\times \Z/4&\Z/2\times \Z/6&(\Z/3\Z)^2&(\Z/4\Z)^2\\
\hline E_G&A_1^{12}&A_1^{14}&A_1^{15}&A_3^4\oplus A_1^4&A_5^3\oplus A_1^3& A_2^8&A_3^6\\
\hline\end{array}\\
\end{array}
\end{eqnarray}

Let $M_G$ be the minimal primitive sublattice of $NS(S_G)$ which
contains the curves $M_i$, then $M_G$ is an overlattice of finite
index $r_G$ of $E_G$ and its properties are the followings \small
$$
\begin{array}{l}
\begin{array}{|c||c|c|c|c|c|c|c|c|c|c|c|c|c|c|}
\hline
G&\Z/2\Z&\Z/3\Z&\Z/4\Z&\Z/5\Z&\Z/6\Z&\Z/7\Z&\Z/8\Z\\
\hline
r_G&2&3&4&5&6&7&8\\
\hline
\rk(M_G)&8&12&14&16&16&18&18\\
\hline
M_G^{\vee}/M_G&(\Z/2\Z)^6&(\Z/3\Z)^4&(\Z/2\Z\times\Z/4\Z)^2&(\Z/5\Z)^2&(\Z/6\Z)^2&(\Z/7\Z)&\Z/4\Z\times \Z/2\Z\\
\hline\end{array}\\
$ $\\

\begin{array}{|c||c|c|c|c|c|c|c|c|c|c|c|c|c|c|}
\hline
G&(\Z/2\Z)^2&(\Z/2\Z)^3&(\Z/2\Z)^4&\Z/2\times \Z/4&\Z/2\times \Z/6&(\Z/3\Z)^2&(\Z/4\Z)^2\\
\hline
r_G&2^2&2^3&2^4&8&12&3^2&4^2\\
\hline
\rk(M_G)&12&14&15&16&18&16&18\\
\hline
M_G^{\vee}/M_G&(\Z/2\Z)^8&(\Z/2\Z)^8&(\Z/2\Z)^7&(\Z/2\Z\times \Z/4\Z)^2&\Z/2\Z\times \Z/6\Z&(\Z/3\Z)^4&(\Z/4\Z)^2\\
\hline\end{array}
\end{array}$$ \normalsize 

The roots of the lattice $M_G$
coincide with the roots of the lattice $E_G$ for all the abelian groups $G$.

By construction $M_G$ is a negative definite lattice primitively embedded in $NS(S_G)$ and
thus $\rho(S_G)\geq 1+\rk(M_G)$.
\end{proposition}
\proof The classification of the Abelian groups acting symplectically on a K3 surface is given in \cite[Theorem 4.5]{Niksympl}, where it is also proved that $S_G$ is a K3 surface if and only if $G$ acts symplectically on $Y_G$. The lattices $E_G$ and $M_G$ are described in \cite[Sections 6 and 7]{Niksympl}. The fact that the root lattices of $M_G$ and of $E_G$ coincide can be checked by a Magma computation as in proof of Proposition \ref{prop: definition lattices FG and KG}. \endproof

We obtain an analogue of Theorem \ref{theorem: X is A/G iff KG is in NS(X)} proved before. 
\begin{theorem}\label{theorem: S is A/G iff MG is in NS(S)} Let $G$ be one of the Abelian groups acting symplectically on a K3 surface. A K3 surface $S_G$ is the desingularization of the quotient $Y_G/G$ for a certain K3 surface $Y_G$ if and only if $M_G$ is primitively embedded in $NS(S_G)$.\end{theorem}
\proof The proof is similar (but easier) to the one of Theorem \ref{theorem: X is A/G iff KG is in NS(X)}. Since the Abelian groups $G$ acting symplectically on a K3 surface are either cyclic or free products of cyclic groups, there is a correspondence between the divisible classes of $M_G$ and covers of $S_G$, given by Section \ref{sec: covers}. So it is immediate to reconstruct the covering surface and its minimal model $Y_G$ from the following data: $S_G$, the lattice $M_G$, the knowledge that certain $(-2)$ classes in $M_G$ represent smooth irreducible rational curves on $S_G$. The latter condition is guaranteed by Lemma \ref{lemma: by Bertin} and the fact that the roots of $E_G$ coincide with the ones of $M_G$, see Proposition \ref{prop: EG and MG}.\endproof

It is not possible to generalize Theorem \ref{theor: it suffices to have the curves on X} or Corollary \ref{cor: it suffices to have an overlattice of  F} to all the Abelian group acting symplectically on a K3 surface. Indeed, for example, there exist K3 surfaces which contain a set of 8 disjoint rational curves, but this set is not divisible by 2, hence these K3 surfaces are not necessarily desingularization of quotient of another K3 surface by $\Z/2\Z$: an example is given by the K3 surface which is the minimal model of the 2:1 cover of $\mathbb{P}^2$ branched along a sextic with 8 nodes. Indeed the cover of $\mathbb{P}^2$ is singular and has 8 singularities of type $A_1$. So on the K3 surface there are 8 disjoint rational curves arising from the desingularization of these singularities. But these curves are not a divisible set: this can be checked considering that the fixed locus of the cover involution is a curve of genus 2 and this determines, by \cite{Nikinvnonsympl}, the N\'eron--Severi group of the K3 surface. It is known that the Theorem \ref{theor: it suffices to have the curves on X} can be extended to the K3 surfaces which contain at least 14 disjoint rational curves, see \cite{GS}.

\begin{rem} The Theorem \ref{theorem: S is A/G iff MG is in NS(S)} was proved for $G=\Z/2\Z$ in \cite[Proposition 2.3]{GSinv} with a different method. The approach used in \cite{GSinv} is strictly based on a careful description of the action induced by a symplectic involution on $\Lambda_{K3}$. This allows one to give stronger results, but a similar description of the action induced by a group of symplectic automorphisms on $\Lambda_{K3}$ is not known for groups $G$ different from $\Z/2\Z$.\end{rem}

The Theorem \ref{theorem: S is A/G iff MG is in NS(S)} allows one to describe the moduli space of the K3 surfaces which are covered by other K3 surfaces in terms of lattice polarized K3 surfaces:

\begin{corollary}\label{corollary: L polarized A/G}
Let $G$ be a finite group acting symplectically on a K3 surface. Let $\mathcal{W}_G$ be the set of lattices $W_G$ satisfying \begin{itemize}\item[(a)] $W_G$ has rank $1+\rk(M_{G})$,\item[(b)] $W_G$ is hyperbolic,\item[(c)] $W_G$ admits a primitive embedding in $\Lambda_{K3}$, \item[(d)] $M_G$ is primitively embedded in $W_G$.\end{itemize} A K3 surface is the desingularization of the quotient of a K3 surface by $G$ if and only if it is a $W_G$-polarized K3 surface for a $W_G\in\mathcal{W}_G$. 

In particular the coarse moduli space of the K3 surface which are desingularization of the quotient $Y/G$ for a K3 surface $Y$ has infinitely many components of dimension $19-\rk(W_G)$.\end{corollary}

In case $G=\Z/2\Z$ all the admissible lattices which appears in $\mathcal{W}_{\Z/2\Z}$ are described in \cite[Proposition 2.1 and Corollary 2.1]{GSinv}. Here we obtain the analogous result for $G=\Z/3\Z$. First we fix the following notation. The lattice $E_{\Z/3\Z}$ is isometric to $A_2^6$. We denote by $M_i^{(j)}$, $i=1,2$ the two curves which generate the $j$-th copy of $A_2$ in $E_G$ and by $d_j:= \left(M_1^{(j)}+2M_2^{(j)}\right)/3$. We can assume that  $M_{\Z/3\Z}$ is generated by the generators of $E_{\Z/3\Z}$ and by the class $\sum_{j=1}^6 d_j$.

\begin{proposition}\label{prop: NS K3 Z3Z}
Let $Y_{\Z/3\Z}$ be a K3 surface which admits a symplectic action of $\Z/3\Z$. Let $S_{\Z/3\Z}$ be the K3 surface desingularization of $\left(Y_{\Z/3\Z}\right)/\left(\Z/3\Z\right)$. Let us assume that $\rho(S_{\Z/3\Z})=13$. There is a primitive embedding of $M_{\Z/3\Z}$ in $NS(S_{\Z/3\Z})$. Let us denote by $H$ a generator of the 1-dimensional subspace of $NS(S_{\Z/3\Z})$ orthogonal to $M_{\Z/3\Z}$ in $NS(S_{\Z/3\Z})$. So $H^2=2d$ for a positive integer $d$ and without loss of generality we can assume that $H$ is pseudoample. Then there are the following possibilities and all of them appear:
\begin{itemize}
\item $d\not\equiv 0\mod 3$: in this case $NS(S_{\Z/3\Z})\simeq \Z H\oplus M_{\Z/3\Z}$;
\item $d\equiv 0\mod 3$: in this case there are two possibilities, either $NS(S_{\Z/3\Z})=\Z H\oplus M_{\Z/3\Z}$ or $NS(S_{\Z/3\Z})$ is an overlattice of index 3 of $\Z H\oplus M_{\Z/3\Z}$. In the latter case $NS(S_{\Z/3\Z})$ is generated by the generators of $M_{\Z/3\Z}$ and by a class $v$. Up to isometries the class $v$ (mod $\Z H\oplus M_{\Z/3\Z}$) is uniquely determined by $d\mod 9$ and it is the following:
\begin{itemize}\item if $d\equiv 0\mod 9$, then $v:=H/3+\sum_{j=1}^3 d_j$;
\item if $d\equiv 3\mod 9$, then $v:= H/3+\sum_{j=1}^2 (d_j)+2 \sum_{h=3}^4(d_h)$;
\item if $d\equiv 6\mod 9$, then $v:=H/3+d_1+2d_2$. \end{itemize}\end{itemize}
\end{proposition}
\proof The proof is based on the lattice theory and is analogous to the one of \cite[Propositions 2.1, 2.2 and Corollary 2.1]{GSinv}.

Let $S_{\Z/3\Z}$ be a K3 surface which is a desingularization of $Y_{\Z/3\Z}/(\Z/3\Z)$ for a certain K3 surface $Y_{\Z/3\Z}$. Then $M_{\Z/3\Z}$ is primitively embedded in $NS(S_{\Z/3\Z})$ and its orthogonal is a positive definite sublattice of rank 1. 

So $NS(S_{\Z/3\Z})$ is an overlattice of finite index, $s$, of $\Z H\oplus M_{\Z/3\Z}$ where $H^2=2d>0$. 
The discriminant group of the lattice $\Z H\oplus M_{\Z/3\Z}$ is $\Z/2d\Z\times (\Z/3\Z)^4$, so the lattice $\Z H\oplus M_{\Z/3\Z}$ has length 5 if $d\equiv 0\mod 3$ and 4 otherwise. A lattice of length at most 5 and of rank 13 admits a primitive embedding in $\Lambda_{K3}$. Thus for each value of $d$ there are K3 surfaces $S_{\Z/3\Z}$ with $NS(S_{\Z/3\Z})\simeq \Z H\oplus M_{\Z/3\Z}$ and so for any value of $d$ there is a K3 surface obtained as quotient of $Y_{\Z/3\Z}$ by $\Z/3\Z$ and such that $NS(S_{\Z/3\Z})\simeq \Z H\oplus M_{\Z/3\Z}$.

Let us now assume that the index $s$ of the inclusion $\Z H\oplus   M_{\Z/3\Z}\hookrightarrow NS(S_{\Z/3\Z})$ is not 1. Then there is a vector $v$ non trivial in $\left(\Z H\oplus   M_{\Z/3\Z}\right)/ NS(S_{\Z/3\Z})$. Since $M_{\Z/3\Z}$ is primitively embedded in $NS(S_{\Z/3\Z})$, the vector $v$ is of the following form: $v:=\frac{1}{s}(H+m)$, where $m\in M_{\Z/3\Z}$ and  $m/s$ is a non trivial element in the discriminant group of $M_{\Z/3\Z}$. This forces $s$ to be 3. The condition $vH=2d/3\in \Z$ forces $d$ to be a multiple of 3. 

In order to identify $v$ we describe the discriminant group of $M_{\Z/3\Z}$. Let us recall that $M_{\Z/3\Z}$ is an overlattice of index 3 of $E_{\Z/3\Z}\simeq A_2^6$. Since the lattice $M_{\Z/3\Z}$ is obtained by $E_{\Z/3\Z}$ adding the vector $\sum_{j=1}^6 d_j$, the vectors in the discriminant group of $M_{\Z/3\Z}$ are the vectors $\sum_{j=1}^6 \alpha_j d_j$ with $\alpha_i\in \Z/3\Z$ such that $\sum_{i=1}^6\alpha_i\equiv 0\mod 3$. 
So the vector $v$ is of the form $H/3+w$ where $w=\sum_{j=1}^6 \alpha_j d_j$ with $\alpha_i\in \Z/3\Z$ such that $\sum_{i=1}^6\alpha_i\equiv 0\mod 3$. The self intersection of $v$ is $2d/9+\sum_{i=1}^6\alpha_i^2(-2/3)$. We observe that $\alpha_i^2$ is either 0, if $\alpha_i$ is 0, or $1$. The number $k:=\sum_{i=1}^6\alpha_i^2$ is the number of $\alpha_i\in\Z/3\Z$ which are different from 0.
The condition $v^2\in 2\Z$ is then equivalent to $2d-6k\equiv 0\mod 18$ and so to $d-3k\equiv 0\mod 9$. Since we already know that $d\equiv 0 \mod 3$, we have that $d$ is equivalent to one of the values $0,3,6\mod 9$. If $d\equiv 0\mod 9$, then $3k\equiv 0\mod 9$, so $k\equiv 0\mod 3$. If $k=0$, then the divisor $H/3$ is contained in $NS(X_{\Z/3\Z})$, which is impossible, since by definition $H$ is a generator of the sublattice of $NS(X_{\Z/3\Z})$, orthogonal to $M_{\Z/3\Z}$. If $k=3$, then, up to a permutation of the indices, the unique choice for $v$ is $v:=H/3+d_1+d_2+d_3$. We observe that in this case the vector $H/3+2(d_1+d_2+d_3)+d_4+d_5+d_6$ is contained in $NS(X_{\Z/3\Z})$, because it is $v+\sum_{i=1}^6d_i$. If $k=6$ a priori we have two possible choices for $v$: either $v:=H/3+d_1+d_2+d_3+d_4+d_5+d_6$ or $H/3+2(d_1+d_2+d_3)+d_4+d_5+d_6$. The first is not admissible, since it implies $H/3\in NS(X_{\Z/3\Z})$. The second one is equivalent to the choice $v:=H/3+d_1+d_2+d_3$. So if $d\equiv 0\mod 9$, then $NS(S_{\Z/3\Z})$ is generated by the generators of $M_{\Z/3\Z}$ and by $v:=H/3+d_1+d_2+d_3$.
Similarly, if $d\equiv 3\mod 9$, then either $k=1$ or $k=4$. Since  $\sum_{i=1}^6\alpha_i\equiv 0\mod 3$, $k=1$ is not admissible, so (up to a permutation of the indices) we can assume that $v:=H/3+d_1+d_2+2d_3+2d_4$. If $d\equiv 6\mod 9$, then either $k=2$ or $k=5$. If $k=2$ we can assume that $v:=H/3+d_1+2d_2$. In this case we observe that the vector $H/3+2d_1+d_3+d_4+d_5+d_6$ is contained in $NS(X_{\Z/3\Z})$, because it is the sum of $v$ and $\sum_{i=1}^6d_i$. But the vector $H/3+2d_1+d_3+d_4+d_5+d_6$ is the unique admissible choice for $v$  (up to a permutation of the indices) with $k=5$. So if $d\equiv 6\mod 9$ we can assume that $v:=H/3+d_1+2d_2$.

\endproof

\begin{rem} There is a clear geometric meaning of $H$ and $d$. Indeed, for every value of $d$ there is a projective model of $S_{\Z/3\Z}$, given by  $\phi_{|H|}:S_{\Z/3\Z}\ra \mathbb{P}(H^0(X,H)^{\vee})$. The image $\phi_{|H|}(S_{\Z/3\Z})$ is a surface with 6 singularities of type $A_2$ and it is in fact the quotient surface $Y_{\Z/3\Z}/\Z/3\Z$. The self intersection of $H$ determines the dimension of the ambient space of $\phi_{|H|}(S_{\Z/3\Z})$, which is $\mathbb{P}^{d+1}$. This is the smallest projective space in which one can describe the quotient $Y_{\Z/3\Z}/\Z/3\Z$.  \end{rem}

{\bf Acknowledgements} I'm grateful to Daniel Huybrechts for several stimulating discussions and for asking me questions which inspired this paper. I thank Bert van Geemen for reading a preliminary version of this paper, for his helpful comments and for many interesting discussions.

\end{document}